\documentclass[10pt]{article}
\usepackage[hmargin=3cm,vmargin=3cm]{geometry}
\usepackage{amsmath,amssymb,paralist,url,hyperref,doi}
\usepackage{palatino,epsfig,latexsym,xfrac} %natbib
\usepackage{setspace}
\usepackage{multirow} 
\usepackage{color}
\usepackage{lineno}
\usepackage{comment}
\usepackage{graphicx}
\usepackage{wrapfig}
\usepackage{balance} 
\def\prob{{\mathbf{Pr}}}
\def\E{\mathop{{\mathbb{E}}}}
\newcommand{\opt}[1]{\textcolor{green}{#1}}
\newcommand{\an}[1]{\textcolor{blue}{#1}}
\newcommand{\td}[1]{\textcolor{red}{#1}}
\usepackage{algorithmic}
\usepackage[linesnumbered,boxed]{algorithm2e} %
\DontPrintSemicolon
\usepackage{authblk}
\title{Addressing Unboundedness in Quadratically-Constrained Mixed-Integer Problems} 
\author{Guy Zepko$^{1,2}$ and Ofer M. Shir$^{2,3}$\footnote{Corresponding author: \href{mailto:ofersh@telhai.ac.il}{\texttt{ofersh@telhai.ac.il}}.}}
\affil{ \begin{small} $^{1}$ Reichman University, Herzliya, Israel\\$^{2}$ Migal Institute, Qiryat Shemona, Israel\\$^{3}$ Tel-Hai College, Upper Galilee, Israel\end{small}}
\date{}
\begin{document}
\maketitle              % typeset the header of the contribution
\begin{abstract}
Mixed-integer (MI) quadratic models subject to quadratic constraints, known as All-Quadratic MI Programs, constitute a challenging class of NP-complete optimization problems. The particular scenario of \textit{unbounded integers} defines a subclass that holds the distinction of being even \emph{undecidable} \cite{Jeroslow1973}. This complexity suggests a possible soft-spot for Mathematical Programming (MP) techniques, which otherwise constitute a good choice to treat MI problems.
We consider the task of minimizing MI convex quadratic objective and constraint functions with unbounded decision variables. 
Given the theoretical weakness of white-box MP solvers to handle such models, we turn to black-box meta-heuristics of the Evolution Strategies (ESs) family, and question their capacity to solve this challenge. Through an empirical assessment of all-quadratic test-cases, across varying Hessian forms and condition numbers, we compare the performance of the CPLEX solver to modern MI ESs, which handle constraints by penalty. 
Our systematic investigation begins where the CPLEX solver encounters difficulties (\textit{timeouts} as the search-space dimensionality increases, $D\gtrsim 30$), and we report in detail on the $D=64$ case.
%on which we report by means of detailed analyses. 
%An interesting duality is evident, when the easiest problem for CPLEX turns to be among the hardest for the ESs, and vice-versa, when an ES repetitively obtains precise optimizers faster than individual CPLEX runs that terminate in a \textit{timeout}.
Overall, the empirical observations confirm that black-box and white-box solvers can be competitive, where CPLEX is evidently outperformed on 13\% of the cases. 
This trend is flipped when unboundedness is amplified by a significant translation of the optima, leading to a totally inferior performance of CPLEX at 83\% of the cases. 
%especially when the \textit{constraint function is separable}, and that two common ESs' mutation operators can effectively handle the integer unboundedness.
We also conclude that conditioning and separability are not intuitive factors in determining the hardness degree of this class of MI problems. 
%where regular versus rough landscape structures can pose mirrored degrees of challenge for MP versus ESs.
\end{abstract}
%
%\textbf{Keywords}: Unbounded integer programs; evolution strategies; integer mutation distributions;\\ CMA-ES with Integer Handling; ILOG-CPLEX; undecidability.
%
%
\section{Introduction}
Global optimization is a fundamental task across a wide range of disciplines, which is concerned with locating the absolute best objective function value within a region that was predefined by constraints.
When the analytical forms of the objective function and constraints are known, the problem is classified as a White-Box Optimization (WBO) problem, which can be solved in a bottom-up manner by utilizing the explicit problem structure and available data. For instance, in numerical optimization, if the continuous objective function is convex and and certain conditions are satisfied, exact solvers can guarantee a solution to the problem \cite{Boyd}.
At the same time, problems with non-convex objective functions have been shown to be NP-hard, making them more difficult to solve \cite{pardalos1991quadratic}.
In contrast, when no information about the objective function is known, the optimization problem is referred to as a Black-Box Optimization (BBO) problem.
In practical optimization problems, a blending of WBO and BBO approaches, referred to as gray-box optimizers \cite{Whitley2017NGG} or hybrid metaheuristics \cite{blum2016hybrid}, are frequently encountered. 
Also, a scheme of BBO problems with so-called \textit{explicit constraints}, where the constraints are unveiled while the objective function remains black-boxed, is another form of blending that has received attention recently \cite{Arnold2023_miConstraints}.   
Notably, the distinction between BBO and WBO is also evident through \textit{model-based optimization}, which targets BBO problems by constructing explicit surrogates and treating them by WBO solvers. A recent study successfully brought this target to fruition on constrained discrete BBO problems \cite{pmlr-v162-papalexopoulos22a}.
WBO is typically approached by formal algorithms and Mathematical Programming (MP) techniques, which are rooted in Theoretical Computer Science \cite{PapadimitriouSteiglitz,moore2011nature}. 
Randomized search heuristics, which are a commonly used approach for solving BBO problems \cite{Wegener2004}, operate by evaluating candidate solutions and using the resulting function values to guide the selection of future search points. They serve as a viable alternative to exact solvers, especially on non-convex landscapes.
Evolution Strategies (ESs; \cite{Baeck2013contemporary,ESchapter2018}) are a family of effective randomized search heuristics for continuous and discrete search-spaces.
\paragraph{Focus}
We focus on Mixed-Integer Optimization (MIO) \cite{floudas1995nonlinear,sahinidis2019mixed,li2013mixed,franco2014mixed}, which involves optimizing an objective function subject to constraints that include both continuous and discrete variables. MIO problems are NP-hard already in the linear case (MILP) and finding an optimal solution is often challenging. 
However, various optimization techniques have been developed to tackle MIO problems, including MP such as branch-\&-bound algorithms, branch-\&-cut algorithms, relaxations to semidefinite programming \cite{Schrijver2011,Boyd}, constraint programming \cite{franco2014mixed}, and meta-heuristics \cite{NACOhandbook,CIhandbook,HeuristicsHandbook}. These techniques can be applied to a wide range of real-world problems, including design optimization, process optimization, hyperparameter optimization, resource allocation, and network design \cite{li2013mixed,floudas1995nonlinear}.
The definition of the feasible integer search-space by, e.g., box-constraints, or otherwise being unbounded, is fundamental in the \textit{complexity} sense. 
A striking theorem, proved by Jeroslow half a century ago \cite{Jeroslow1973}, states that \textit{unbounded quadratically-constrained integer programs are undecidable}. 
Accordingly, WBO solvers are challenged in practice when handling this class of problems, due to the nature of their operation (i.e., the mathematical barriers in exploiting the explicit model when subject to integer unboundedness). 
Meta-heuristics from the BBO perspective obviously face the same computational complexity when treating this class of problems, but their operation is not necessarily affected by this unbounded integer search. It has, nevertheless, implications on the choice of techniques. In particular, it eliminates the possibility to employ standard Genetic Algorithms \cite{Goldberg}, which require to encode the feasible space by means of fixed-size bitstrings, or any other meta-heuristic that requires explicit boundaries for its operation (e.g., the so-called CMA-ES with ``Margins'' (\texttt{cma-wM}; \cite{CMAESwM2022}). The careful reader should distinguish the naming 'margins', which is rooted in probability treatment, from the practical issue of integer boundaries, whose explicit definition is a prerequisite for the operation of this technique).
\paragraph{Contribution}
MI ESs (MIESs) were already investigated three decades ago by \cite{rudolph1994evolutionary} and \cite{Baeck95MIES}, and yet, the research outcome of this MIES thread has not widely developed for most years, with a limited number of published studies \cite{RUILI_PPSN2008,Reehuis2010,li2013mixed}, while enjoying a growing interest only recently \cite{miBBOB2019,CMAESwM2022,cmaIHgecco2023,Arnold2023_miConstraints}. 
Importantly, certain fundamental questions still remain unanswered, despite the evident \textit{success in practice} of modern ESs to solve MI problems (e.g., \texttt{cma-wM} \cite{CMAESwM2022} for bounded spaces, or another CMA-ES variant, which was adapted to handle integers \cite{cmaIHgecco2023} (\texttt{cma-IH}), for the unbounded case). 
An open question concerns the effectiveness of employing the \textit{normal} distribution for integer mutations, in either bounded or unbounded spaces.
%Essentially, Rudolph proved the rationale of using the \textit{geometric} distribution for mutations in unbounded integer optimization \cite{rudolph1994evolutionary}, but somehow the normal distribution remains the omni-hammer in the ES community, even when turning to discrete search-spaces. 
The contribution of the current study is twofold -- (i) assessing whether the theoretical undecidability of unbounded quadratically-constrained MI problems constitutes a WBO weakness, and thus an opportunity for BBO in practice, and (ii) assessing the effectiveness of the ES mutation distribution when handling unbounded integers.

\paragraph{Paper Organization} %\td{TODO}
Next, Section \ref{sec:formulation} will formally state the problem, distill the research questions and mention the concrete aims.
Section \ref{sec:approach} will then present the utilized techniques abd specify the experimental setup. 
The numerical results will be presented in Section \ref{sec:results}, and finally, Section \ref{sec:summary} will conclude and summarize this study.
\section{Problem Statement and Formulation}\label{sec:formulation}
A standard quadratic program (QP) is formally defined as the minimization of a quadratic objective function whose decision variables are placed within the unit simplex ($\Delta_D := \left\{\vec{x}\in \mathbb{R}_{+}^D : \vec{e}^T \vec{x}=1 \right\}$, with the vector $\vec{e}\in \mathbb{R}^D$ of all ones):
\begin{equation}\label{eq:standardQP}
     \displaystyle \textrm{minimize}_{\vec{x} \in \Delta_D} \quad \vec{x}^T \mathbf{H} \vec{x}, 
\end{equation}
and where $\mathbf{H}$ is a symmetric $D\times D$ real matrix. Non-homogeneous quadratic forms are easily treated by a trivial reformulation \cite{GY2021}.  
The introduction of integer decision variables renders the problem MI, standardly denoted as MIQP.
Solving Eq.\ \ref{eq:standardQP} is already NP-complete \cite{QPisNP_1990} -- i.e., this so-called ``pure-QP'' of continuous variables is a hard problem even before introducing the integer decision variables.
When such integers are introduced, the resultant MIQP is clearly a harder computational problem, but progress in addressing it in practice has been achieved \cite{MIQP_CPLEX}. 
One of the important aspects of the MIQP branch is the fact that it constitutes the next-step from the well-established MILP branch (see, e.g., \cite{BixbyCPLEX}) toward the generalized Mixed-Integer Nonlinear Programming (MINLP) branch \cite{belotti_kirches_leyffer_linderoth_luedtke_mahajan_2013}.
Finally, another degree of complexity is introduced to the model when the formulation encompasses quadratic constraint terms, and it is then called Quadratically-Constrained QP (QCQP) (also known as \textit{all-quadratic programs} \cite{Raber1998,ZHAO2017159}). In the current context of MI optimization, we will denote such problems as MIQCQP.
Quadratic models, either pure-QP or MIQP, arise in a large variety of problems, ranging from Portfolio optimization (Markowitz's original formulation \cite{Markowitz52} and its extensions \cite{PortfolioMINLP}) and resource allocation \cite{ResAllocProb88}, to population genetics \cite{Kingman1961} and game theory \cite{Bomze2002}.

\subsubsection{WBO: MIQCQP as a Soft-Spot}
Jeroslow proved the undecidability of unbounded integer programs with quadratic constraints \cite{Jeroslow1973}.
Despite this proven problem hardness, there has been much practical progress in treating MIQP in general \cite{MIQP_CPLEX} and MIQCQP in particular (e.g., by reformulation to a bilinear programming problem with integer variables \cite{ZHAO2017159}, or by diverting to Mixed-Integer Second-Order Cone Programming \cite{MISOCP2013} when the model permits).
%\footnote{For pragmatic purposes, convenient implementation interfaces are available -- see, e.g., \cite{CPLEX-MIQCP}.}
Although many linearization and reformulation techniques exist \cite{Glover2004Comparisons,Linearization2022}, a generalized unbounded MIQP cannot be linearized.\footnote{Even if the integers $x_i~:~i\in I$ may be linearized using auxiliary binaries \cite{BIL2008}, the multiplication of two \textit{unbounded} decision variables within $\vec{x}$ cannot be linearized.} 
%(but could be piecewise approximated upon bounding and separation \cite{Linearization2022}).}
WBO solvers are typically challenged when treating MIQCQP, an angle that we intend to explore.

\subsubsection{BBO: On Integer Mutations and Distributions}
%Studying the effectiveness of discrete randomized heuristic search in integer search-spaces dates back to the 1960's and 1970's. 
%Kelahan and Gaddy proposed to use a bilateral power distribution that adaptively shrank by a geometric schedule \cite{Kelahan1978-RandomSearchMIO}.
Rudolph questioned the suitability of random distributions for evolutionary mutations in unbounded integer spaces \cite{rudolph1994evolutionary}. In particular, Rudolph showed that integer mutations utilizing the geometric distribution possess maximum entropy, and demonstrated the effectiveness of such mutations in practice.
That study specifically proposed to use the doubly-geometric mutation, which possesses a symmetric distribution with respect to 0. 
It is achieved by drawing two random variables, $\left\{ g_1,g_2\right\}$, according to the geometric distribution $(\jmath=1,2)$: $\prob\left\{ g_{\jmath} = k \right\} = p\cdot \left(1-p \right)^k$ , and taking their difference, $z=g_1-g_2$.
%The probability function of $z$ reads,
%$$ \prob\left\{ z = k \right\} = \frac{p}{2-p} \cdot \left(1-p \right)^{\left|k\right|},\quad k \in \mathbb{Z} $$
%with $\E \left[z\right]=0$ and $\textrm{VAR}[z]=2(1-p)/p^2$. 
%Next, to generalize to the multivariate case of an $n_z$-dimensional mutation vector $\vec{z}$, Rudolph showed that by taking $n_z$ stochastically independent random variables, following the doubly-geometric distribution, the properties of symmetry and maximal entropy are kept. 
%and the randomly generated variable is calculated via  $\mathcal{G}_{z}\left(0,p\right):=  g_1 - g_2$.\\ 
Importantly, when generalized to the multivariate case of an $n_z$-dimensional mutation vector $\vec{z}$, each dimension is drawn individually as such, 
$\mathcal{G}_{n_z}:g_1 - g_2$, 
yet the distribution as a whole could be \textit{controlled by the mean step-size}, $S=\E \left[ \|\vec{z} \|_1\right]$:\\
\begin{equation} 
\displaystyle p = 1 - \frac{S/n_z}{\sqrt{\left(1+\left(S/n_z \right)^2 \right)} + 1} \quad \Longleftrightarrow \quad S = n_z\cdot \frac{2(1-p)}{p(2-p)}.
\end{equation}
In practice, each random variable is drawn by the following calculation ($g_{\jmath}$ are the geometrically distributed random variables, both with parameter $p$):
\begin{equation}\label{eq:doubleGeometric}
     \begin{array}{l}
         \displaystyle g_{\jmath} \longleftarrow \left\lfloor \frac{\textrm{log} \left(1- \mathcal{U}\left(0,1\right)\right)}{\textrm{log} \left(1- p\right)}\right\rfloor \quad \jmath=1,2 \\
         \displaystyle \mathcal{G}_{n_z}\left(0,p\right):=  g_1 - g_2. 
     \end{array}
\end{equation}
The mentioned geometric distribution underlies classical MIESs \cite{Baeck95MIES,RUILI_PPSN2008}, whereas modern ESs, such as the renowned CMA-ES, were successfully modified to handle integer mutations while maintaining adherence to the normal distribution. 
We are interested in assessing the effectiveness of these two mutation distributions when handling unbounded integer search.

\subsubsection{Research Questions and Concrete Aims}%\td{TODO: research question}
We present our research questions:
\begin{quote}
    {\it Do unbounded MIQCQP models constitute a weakness of WBO solvers in practice? How do MIESs compare, as BBO meta-heuristics representatives, in treating such models across different shapes of quadratic forms? %And to what extent does the underlying mutation distribution influence their effectiveness?
 }
\end{quote}% 
%\subsection{Summary: Concrete Aims}
Overall, we plan to investigate the following family of MI all-quadratic problems ($E$ serves as the parametric constraint level at the quadratic inequality constraints): %($n_r$ real-valued and $n_z$ integer decision variables; $D:=n_r+n_z$): %%($k\in \left\{0,1\right\}$:%; the separation of the search-sub-spaces, $\vec{x}\in\mathbb{R}^{n_r}, \vec{z}\in \mathbb{Z}^{n_z}$ is deliberately explicit):
\begin{align}\label{eq:ourQP}
\boxed{
\begin{array}{ll}
\medskip 
\displaystyle \textrm{minimize}_{\vec{x}} & f\left(\vec{x}\right):= (\vec{x}-\vec{\xi}_{0})^T \cdot \mathbf{H}_{f} \cdot (\vec{x}-\vec{\xi}_{0})\\
\medskip
\displaystyle \textrm{subject to:}& g\left(\vec{x}\right):= (\vec{x}-\vec{\xi}_{1})^T \cdot \mathbf{H}_{g} \cdot (\vec{x}-\vec{\xi}_{1}) \leq E\\
%\displaystyle & \vec{x}^T \vec{x} \leq K \\
\displaystyle & \vec{x}\in\mathbb{R}^{D}\\ 
\displaystyle & x_i\in \mathbb{Z} \quad \forall i \in I,
\end{array} }
\end{align}
where the $D$-dimensional decision vector $\vec{x}$ is constructed by $n_r$ real-valued decision variables followed by $n_z$ integer decision variables that are defined by the so-called \textit{index set} $I:=\left\{n_r+1,\ldots,n_r+n_z \right\}$: $\forall i \in I \quad x_i \in \mathbb{Z}$.
\section{Approach and Methodology}\label{sec:approach}
We would like to conduct a systematic comparison between an exact WBO solver to representative ESs in a \textit{fixed-budget} (resources) plan.
To this end, we consider the commercial CPLEX solver as an MP exact method (and particularly as an MIQCQP solver), versus two MIESs. 
We begin this section by shortly presenting the techniques, and then describing our numerical setup.
\subsection{Solvers}
\subsubsection*{CPLEX}%{CPLEX as an MIQCQP solver}
IBM ILOG CPLEX constitutes a broad MP environment, with a large variety of state-of-the-art algorithms under the hood. Given linear or quadratic optimization models, either pure, mixed-integer or all-integer, it utilizes an automated algorithms selection procedure that guides the engine to employ the most suitable sub-routines (e.g., start solving a pure-LP with the dual-Simplex algorithm, and shift to other techniques when certain conditions arise). 
In our case (MIQCQP with convex objective and constraint functions \cite{CPLEX-MIQCP}), the engine's default selection is a Branch-\&-Cut scheme relying on a QP solver \cite{best2017quadratic}.
\subsubsection*{\texttt{mies}}
This ES was originally defined to treat altogether real-valued, integer ordinal, and categorical decision variables. 
We discard herein the treatment of categorical variables. 
Like most ESs, its operation is well defined by its self-adaptive mutation operator (Algorithm \ref{algo:mutate}): $\left\{\vec{x},\vec{z}\right\}$ are the real-valued and integer decision vectors, respectively, and $\left\{\vec{s},\vec{q}\right\}$ are their strategy parameters, respectively. 
The normal distribution plays the dominant role of the real-valued update steps, whereas the integer decision variables are mutated by adding \textit{doubly geometrically} distributed random numbers, following Eq.~\ref{eq:doubleGeometric}.
\IncMargin{0.25em}
{\LinesNotNumbered
\begin{algorithm}[h]
\caption{The self-adaptive mutation operator utilized by the \texttt{mies}: 
$\left\{\vec{x},\vec{s}\right\}$ are the real-valued decision variables and strategy parameters, respectively.
$\left\{\vec{z},\vec{q}\right\}$ are the integer decision variables and strategy parameters, respectively. 
$\mathcal{N}$ and $\mathcal{G}$ denote the normal and the geometric distributions, respectively (for the latter see \eqref{eq:doubleGeometric}). 
We set $\varepsilon:=10^{-5}$ and do not enforce boundary constraints. \label{algo:mutate}}
\normalsize 
\Indm\SetKwFunction{Func}{\textbf{mies::mutate}}
\Func{$\vec{x},~\vec{s},~n_r,~\vec{z},~\vec{q},~n_z$}\;
\Indp
  \quad \quad \textbf{/* real-valued decision variables */} \;
  $\mathcal{N}_g^{(r)}\leftarrow \mathcal{N}\left(0,1\right),~\tau_g^{(r)}\leftarrow\frac{1}{\sqrt{2\cdot n_r}},~\tau_{\ell}^{(r)}\leftarrow\frac{1}{\sqrt{2\cdot \sqrt{n_r}}}$ \;
  \For{$i=1,\ldots,n_r$} { 
        $s_i^{\prime}\longleftarrow \max\left(\varepsilon,~s_i\cdot \exp\left\{\tau_g^{(r)}\cdot\mathcal{N}_g^{(r)} + \tau_{\ell}^{(r)}\cdot\mathcal{N}\left(0,1\right)\right\} \right)$ \;
        $x_i^{\prime}\longleftarrow x_i + \mathcal{N}\left(0,s_i^{\prime}\right)$ \;
    }
    \quad \quad \textbf{/* integer decision variables */} \;
    $\mathcal{N}_g^{(z)}\leftarrow \mathcal{N}\left(0,1\right),~\tau_g^{(z)}\leftarrow\frac{1}{\sqrt{2\cdot n_z}},~\tau_{\ell}^{(z)}\leftarrow\frac{1}{\sqrt{2\cdot \sqrt{n_z}}}$ \;
  \For{$i=1,\ldots,n_z$} { 
        $q_i^{\prime}\longleftarrow \max\left( 1, ~q_i\cdot \exp\left\{\tau_g^{(z)}\cdot\mathcal{N}_g^{(z)} + \tau_{\ell}^{(z)}\cdot\mathcal{N}\left(0,1\right)\right\} \right)$ \;
        %$\psi \longleftarrow 1 - \left(q_i^{\prime}/n_z\right) \cdot \left(1 + \sqrt{1 + \left(\frac{q_i^{\prime}}{n_z}\right)^2 } \right)^{-1}$ \;
        %$g_1 \longleftarrow \lfloor \frac{\textrm{ln} \left(1- \mathcal{U}\left(0,1\right)\right)}{\textrm{ln} \left(1- \psi\right)}\rfloor$ ~ $g_2 \longleftarrow \lfloor \frac{\textrm{ln} \left(1- \mathcal{U}\left(0,1\right)\right)}{\textrm{ln} \left(1- \psi\right)}\rfloor$ \;
        %$z_i^{\prime}\longleftarrow z_i + g_1 - g_2$ \;
        $z_i^{\prime}\longleftarrow z_i + \mathcal{G}_{n_z}\left(0,q_i^{\prime}\right)$ // see Eq.~\ref{eq:doubleGeometric} \;
    }
\Return {$\left\{\vec{x}^{\prime},~\vec{s}^{\prime},~\vec{z}^{\prime},\vec{q}^{\prime} \right\}$}
\end{algorithm}}
\DecMargin{0.25em}
\subsubsection*{\texttt{cma-IH}} %\td{TODO: elaborate on the CMA-ES: Margins versus Integer Handling}
The CMA-ES is a modern ES \cite{Baeck2013contemporary} %, which was released as the fourth generation of derandomization ESs \cite{Hansen01completely}, and 
that has enjoyed a broad success in global optimization of continuous problems. 
Its operation is driven by two mechanisms that statistically learn past mutations: updating the covariance matrix $\mathbf{C}$, which is central to landscape maneuvering, and controlling the step-size $\sigma$. 
Those mechanisms were originally defined for continuous landscapes, and experienced performance issues when deployed on MI problems.
%and then were iteratively refined throughout the years in light of a sheer volume of research \cite{BBOB_GECCO2010}.
Most dominantly, the major issue was identified as \textit{stagnation} whenever the discrete decision variables get stuck on the landscape's integer plateaus, and affect the heuristic's progress. To remedy this malfunction on MI landscapes, the notion of margins was introduced to fix the probability of mutating to another integer, resulting in the so-called \texttt{cma-wM} \cite{CMAESwM2022}. However, this fix requires the explicit mapping of the bounded integer space, and thus becomes irrelevant in the unbounded case. A more simplistic fix for handling integers, which we denote as \texttt{cma-IH} and adopt for our study on unboundedness, treats the stagnation by setting a lower bound on the mutation variances \cite{cmaIHgecco2023}. Importantly, the \texttt{cma-IH} handles the entire set of decision variables by means of a unified covariance matrix, which facilitates the $D$-dimensional normally-distributed mutations $\vec{z}_k$, and then applies \textit{integer rounding} to the variables belonging to the index set $I$:\\ %The mutation step ($k=1,\ldots,\lambda$) reads
%\begin{align*}
$\displaystyle \vec{x}^{\prime}_k \sim \mathcal{N}(\vec{m},\sigma^2 \mathbf{C}) = 
\vec{m} + \sigma \cdot \mathcal{N}(\vec{0},\mathbf{C}) = \vec{m} + \sigma \mathbf{C}^{\frac{1}{2}} \vec{z}_k, \quad k=1,\ldots,\lambda$ \\ 
$\displaystyle x^{\prime}_{k,i} \longleftarrow \texttt{round}\left(x^{\prime}_{k,i} \right) \quad \forall i \in I$.
%\end{align*}
%with $\vec{z}_k \sim \vec{\mathcal{N}}\left(0,1\right)$ i.i.d.\ for each individual $k$.

\subsection{Experimental Setup}\label{sec:setup}
To address our research questions, we define MIQCQP models on which we will assess the behavior of representative WBO and BBO techniques by means of numerical simulations. We elaborate on the choices made in our setup.
\subsubsection{Landscape Choice and Instance Generation}
To set up a testing framework, we choose concrete quadratic landscapes, and justify in what follows the selection process. 
We consider one separable and two non-separable Hessian matrices, %($n \in \{n_r,~n_z\}$), 
which play roles as both the objective and the constraint functions:
\begin{enumerate}[(H-1)]
\item Cigar:~$\left( \mathcal{H}_{\textrm{cigar}}\right)_{11} = 1,~\left( \mathcal{H}_{\textrm{cigar}}\right)_{ii} = c~~~i=2,\ldots,n$
\item Rotated Ellipse:~$\mathcal{H}_{\textrm{RE}} = \mathcal{R}\mathcal{H}_{\textrm{ellipse}} \mathcal{R}^{-1}$, where $\mathcal{R}$ is the rotation by $\approx \frac{\pi}{4}$ radians in the plane spanned by $(1,0,1,0,\ldots)^T$ and $(0,1,0,1,\ldots)^T$;
\item Hadamard Ellipse: $\mathcal{H}_{\textrm{HE}} = \mathcal{S}\mathcal{H}_{\textrm{ellipse}} \mathcal{S}^{-1}$, where the rotation constitutes the normalized Hadamard matrix, $\mathcal{S}:=\textrm{Hadamard}(n)/\sqrt{n}$.
\end{enumerate}
$\mathcal{H}_{\textrm{ellipse}}$ is the separable ellipse, with $\left( \mathcal{H}_{\textrm{ellipse}}\right)_{ii} = c^{\frac{i-1}{n-1}}~~~i=1,\ldots,n$ and zero elsewhere. $c$ denotes a \textit{parametric condition number} throughout the test-cases.
This choice is justified by the evident challenge that these Hessians introduced to both approaches in practice: 
Selecting (H-1) was based on preliminary runs with the CPLEX solver which indicated \textit{timeouts} and thus reflected a WBO challenge. 
Equivalently, (H-2) and (H-3) both introduced major challenges to the MIESs. 
Altogether, we investigate the following nine test-cases, where $\mathbf{H}_f$ refers to the Hessian of the quadratic objective function, and $\mathbf{H}_g$ refers to the Hessian of the quadratic constraint function \eqref{eq:ourQP}:\\
\begin{center}
 \boxed{\begin{tabular}{l| l l l }
\multicolumn{1}{l}{} & \multicolumn{3}{c}{$\textbf{H}_g$} \\ \cline{2-4}
$\textbf{H}_f$ & \text{Cigar} & \text{RotEllipse}  & \text{HadEllipse} \\ \hline
\text{Cigar} & TC-0: $\left(\mathcal{H}_{\text{cigar}},\mathcal{H}_{\text{cigar}}\right)$ & TC-1: $\left(\mathcal{H}_{\text{cigar}},\mathcal{H}_{\text{RE}}\right)$  & TC-2: $\left(\mathcal{H}_{\text{cigar}},\mathcal{H}_{\text{HE}}\right)$\\
\text{RotEllipse} & TC-3: $\left(\mathcal{H}_{\text{RE}},\mathcal{H}_{\text{cigar}}\right)$ & TC-4: $\left(\mathcal{H}_{\text{RE}},\mathcal{H}_{\text{RE}}\right)$  & TC-5: $\left(\mathcal{H}_{\text{RE}},\mathcal{H}_{\text{HE}}\right)$\\
\text{HadEllipse} & TC-6: $\left(\mathcal{H}_{\text{HE}},\mathcal{H}_{\text{cigar}}\right)$ & TC-7: $\left(\mathcal{H}_{\text{HE}},\mathcal{H}_{\text{RE}}\right)$  & TC-8: $\left(\mathcal{H}_{\text{HE}},\mathcal{H}_{\text{HE}}\right)$\\
\end{tabular}}
\end{center}

\begin{figure}
    \centering
%    \medskip
    \begin{tabular}{ c | c | c}
    \hline
    \includegraphics[width=0.3\columnwidth]{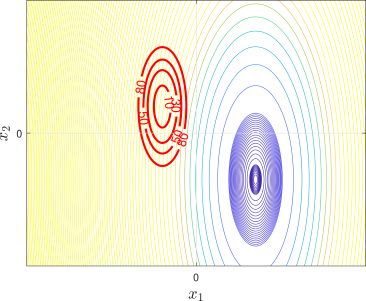} & \includegraphics[width=0.3\columnwidth]{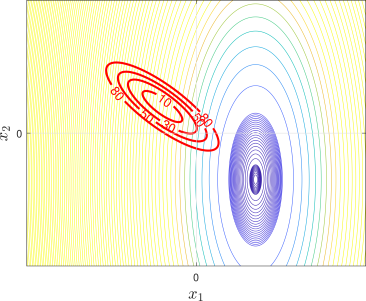} & \includegraphics[width=0.3\columnwidth]{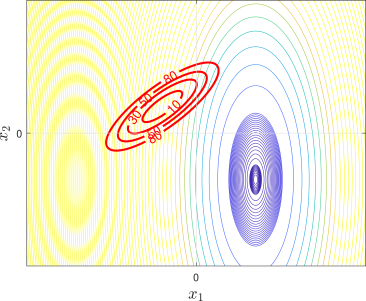}\\
    \hline
    \includegraphics[width=0.3\columnwidth]{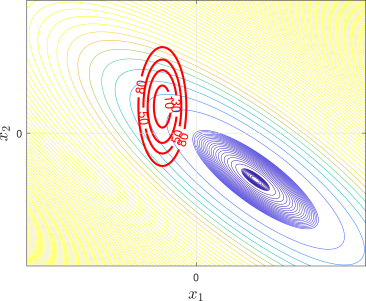} & \includegraphics[width=0.3\columnwidth]{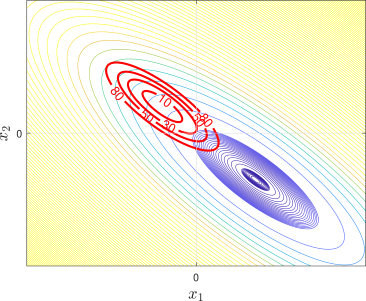} & \includegraphics[width=0.3\columnwidth]{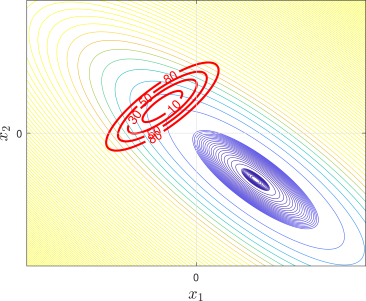}\\
    \hline
    \includegraphics[width=0.3\columnwidth]{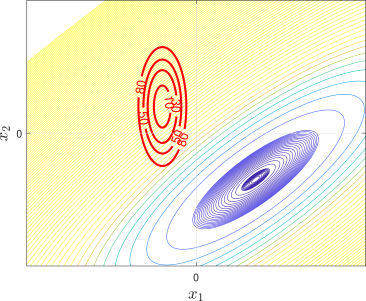} & \includegraphics[width=0.3\columnwidth]{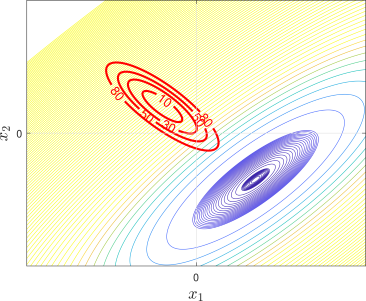} & \includegraphics[width=0.3\columnwidth]{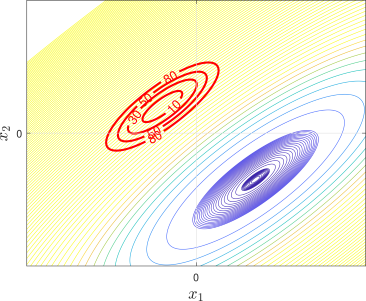}\\
    \hline
    \end{tabular}
    \caption{Summary of the $3\times 3$ test-cases to be investigated as ordered pairs of $\left(\mathbf{H}_f,\mathbf{H}_g\right)$: %[TOP] an index table for the nine test-cases; [BOTTOM] 
    Contour plots of the corresponding \textbf{relaxed} 2D problem instances, depicting the underlying objective function (as the colorful landscape) and the \textbf{quadratic constraint function} at four \textit{level-sets} (thick red curves at $E=\left\{10,30,50,80 \right\}$ \eqref{eq:ourQP}), which defines the feasible ellipsoids on each landscape. The conditioning is set to $c=10$. Notably, the RotEllipse+RotEllipse (TC-4) constitutes a pair of aligned ellipsoids along their long axes, which remain aligned with the increasing conditioning. The Hadamard rotation generates a different impact -- the HadEllipse+HadEllipse (TC-8) forms a pair of parallel ellipsoids that diverge,
 creating a ridge-like effect, as the conditioning increases.}\label{fig:landscapes-gallery}
\end{figure}

\subsubsection{Concrete Calculations}
We consider an equal portion of continuous versus discrete decision variables, $D:=n_r+n_z=2\cdot n$.
We set the following points in $\mathbb{R}^{D}$ about which the quadratic models are centered:
\begin{align}\label{eq:xi}
\displaystyle \vec{\xi}_0:= & \left(+7,-7,+7,-7,...,+7,-7\right)^T &
\displaystyle \vec{\xi}_1:= & \left(-4,+4,-4,+4,...,-4,+4\right)^T.
\end{align}
For the separable Cigar, we construct the $D\times D$-dimensional Hessian matrix as a concatenation of two $n\times n$-dimensional matrices, in order to introduce the conditioning effect to both types of variables:
$$\mathbf{H}:= \begin{pmatrix} \mathcal{H} \quad \mathbf{0} \\ ~ \\ \mathbf{0} \quad \mathcal{H} \end{pmatrix}.$$
Once the $D\times D$-dimensional Hessian matrix is set per the objective function $f$ (denoted as $\mathbf{H}_f$) as well as per the constraint function $g$ (denoted as $\mathbf{H}_g$), their values are calculated as follows ($\vec{y}_k$ denote the translated $\vec{x}$ with respect to $\vec{\xi}_k$):
\begin{align}
\displaystyle f =  \frac{1}{c} \cdot \vec{y}_0^T \mathbf{H}_f \vec{y}_0, \quad & \displaystyle \quad g =  \frac{1}{c} \cdot \vec{y}_1^T \mathbf{H}_g \vec{y}_1.
\end{align}
Importantly, the MIESs handle the inequality constraint of $g$ by means of a \textit{penalty term}, and practically treat the following \textit{cost function}:
\begin{align}\label{eq:cost}
\textbf{\textrm{[BBO cost function:]}} \quad & \quad f+ 10^4 D^2 \cdot \Theta(g-E)\cdot (g-E)^2 \mapsto \min,    
\end{align}
where $\Theta$ denotes the Heaviside step function.

\subsubsection{Preliminary: Mixed-Integer Sphere-constrained-by-Sphere}
We conducted preliminary runs to simulate WBO problem-solving by CPLEX of the simplest MIQCQP instance, subject to various parabolic levels, $E=\{10,30,50,80\}$, across increasing dimensionalities, $D=\{4,8,16,32,64\}$, each limited to 1 hour: % i.e., the program of Eq. \ref{eq:ourQP} with 
\begin{align}\label{eq:SphereBySphere}
\displaystyle \textrm{minimize}_{\vec{x}} ~~  \vec{y}_0^T \cdot \vec{y}_0 \quad & \displaystyle \quad \displaystyle \textrm{subject to: } \vec{y}_1^T \cdot \vec{y}_1 \leq E.
\end{align}
We summarize the outcome of these preliminary runs in Table \ref{tab:preliminaryCPLEXsphere}. 
Evidently, the CPLEX solver is challenged on this MIQCQP test-case as of dimension $D=32$, when it occasionally terminates upon exceeding the time-limit. We therefore set up our experiments to consider this dimension as a starting point.
\begin{table}[]
    \centering
    \caption{CPLEX preliminary runs on the Sphere-constrained-by-Sphere case (Eq. \ref{eq:SphereBySphere}).
    The entries specify the optimization outcome within a time-limit of 1 computation hour -- \texttt{optimal}: an optimal solution was attained; \texttt{tolerance}: a solution was located within the relative gap tolerance; and \texttt{time-out} when exceeded the time-limit without indication regarding optimality of the candidate solution.}
    \label{tab:preliminaryCPLEXsphere}
    \begin{small}
    \begin{tabular}{l|c|c|c|c}
    \hline
         Dimen. & $E=10$ & $E=30$ & $E=50$ & $E=80$ \\
         \hline
         $D=4$ & \opt{\texttt{optimal}} &  \opt{\texttt{optimal}} &  \opt{\texttt{optimal}} &  \opt{\texttt{optimal}}\\
         $D=8$ &  \opt{\texttt{optimal}} & \an{\texttt{tolerance}} & \an{\texttt{tolerance}} &  \opt{\texttt{optimal}}\\
         $D=16$ & \an{\texttt{tolerance}} & \an{\texttt{tolerance}} & \an{\texttt{tolerance}} & \an{\texttt{tolerance}}\\
         $D=32$ & \td{\texttt{time-out}} & \an{\texttt{tolerance}} &\an{\texttt{tolerance}} & \td{\texttt{time-out}}\\
         $D=64$ & \td{\texttt{time-out}} & \td{\texttt{time-out}} & \an{\texttt{tolerance}} & \an{\texttt{tolerance}}\\
         \hline
    \end{tabular}
    \end{small}
\end{table}
\subsubsection{Problem Instances} %Generation
To generate concrete problem instances (i.e., instantiate $k$ and $\ell$ in Eq.\ \ref{eq:ourQP}), and given the preliminary runs, we decided to explore three dimensions, $D=\{32,64,128\}$, with constraints level, $E=\{10,30,50,80\}$, across six conditioning, $c=\{10,100,\ldots,10^6\}$. 
Considering the double-Hessian combinations, this setup yields altogether 216 problem instances per dimension.
\textbf{Importantly, we will focus in Section \ref{sec:results} on experimenting in depth the 64-dimensional use-case (i.e., $n=32,~D=64$). }
Preliminary runs reflected equivalent algorithmic behavior on the various dimensions, and yet, this particular dimension is selected for being an interesting tradeoff between high dimensionality to known scalability issues of ESs (see, e.g., \cite{Shir-ML-JHEUR2021}).
%Preliminary runs reflected equivalent algorithmic behavior on the various dimensions, and yet, this particular dimension is selected for being an interesting tradeoff between high dimensionality to known scalability issues of ESs \cite{Shir-ML-JHEUR2021}.
%
\subsubsection{Numerical Setup}
%\td{TODO: adjust} 
%Preliminary runs indicated that the CPLEX required long computation times for obtaining an exact solution (``optimal''), or so-called ``integer optimal up to the tolerance''. Thus, a pragmatic time limit was much needed.
%It may be argued that the CPLEX benefits from a CPU advantage over the MIESs, which require repetitions due to their stochasticity.
Adhering to a \textit{fixed-budget} plan, we designate similar computational resources for each method: 1 hour per problem instance at $D=\{32,64\}$, and 2 hours at $D=128$.
However, due to the ESs' stochasticity, we allow a single CPLEX run versus 10 MIESs' stochastic runs (each granted a budget of $10^6$ function evaluations, which was doubled at $D=128$). \\
Next, we provide the technical specifications of our numerical simulations.\footnote{The source code will be provided upon request.}
  
\paragraph{WBO Setup} IBM ILOG CPLEX Optimization Studio 12.8 facilitated the WBO problem solving, modeled and run in OPL.
    All the experiments were run using the \texttt{Python API} (for sequential execution) and executed on Windows Intel(R) Xeon(R) CPU E5-1620 v4 @ 3.50GHz with 16 processing units. The relative MIP optimality gap was set to $10^{-3}$ (\texttt{cplex.epgap = 0.001}), and the \texttt{polish} procedure \cite{ILOGpolishing} was enabled after reaching an integer solution. 
    
\paragraph{BBO Setup} 
    The MIESs were implemented in \texttt{python3}; our \texttt{mies} implementation and parameter settings typically followed \cite{Reehuis2010}, whereas the \texttt{cma-IH} was deployed using the \texttt{pycma} package (version \texttt{3.3.0}) with the default integer handling option. Preliminary runs indicated that the recombination operator was \textit{disruptive} for the \texttt{mies}, so it was disabled.
    All runs were deployed on Linux Intel(R) Xeon(R) CPU E5-4669 v4 @ 3.0GHz with 88 processing units.

\section{Numerical Results}\label{sec:results}
Performance comparisons between an exact solver to a randomized search heuristic can be evaluated differently, for instance by considering either the best attained result of the heuristic or the median value of its runs. Generally speaking, questioning the fairness of such performance evaluations may be open for a debate.
In the current study we designed the experimental campaign in such a manner that the computational resources of a single CPLEX run were equivalent in practice to 10 parallel runs of the MIESs, which justifies the empirical comparison of the \textbf{best} attained MIESs' values. 
At the same time, we will also explore the attained \textbf{median} values to assess the MIESs' behavior and possibly infer on the search challenge. 
\subsection{Benchmark on MIQCQP Test Problems}\label{sec:benchmark}
We report on the numerical outcome of the simulations prescribed in Section \ref{sec:setup} and analyze them.
We focus on presenting the results per $D=64$. 
%, and explicitly mention observations per $D=32$ and $D=128$ only when they exhibited different trends.\footnote{The entire dataset of raw results and plots will be published upon acceptance.}
In terms of performance, it is important to note that the CPLEX always terminates due to \td{\texttt{time-out}} throughout the test-cases, except for the (TC-4) problem (RotEllipse+RotEllipse) with high conditioning instances, on which it terminated instantly with a status \an{\texttt{tolerance}}. %Figures \ref{fig:D64sep} and \ref{fig:D64nonsep} provide  %
We first report on high-level performance with respect to attained best values (\ref{sec:bestRes}), and then we delve into the statistics of the runs (\ref{sec:medRes}).
\begin{figure}
    \centering
    \includegraphics[width=0.48\columnwidth]{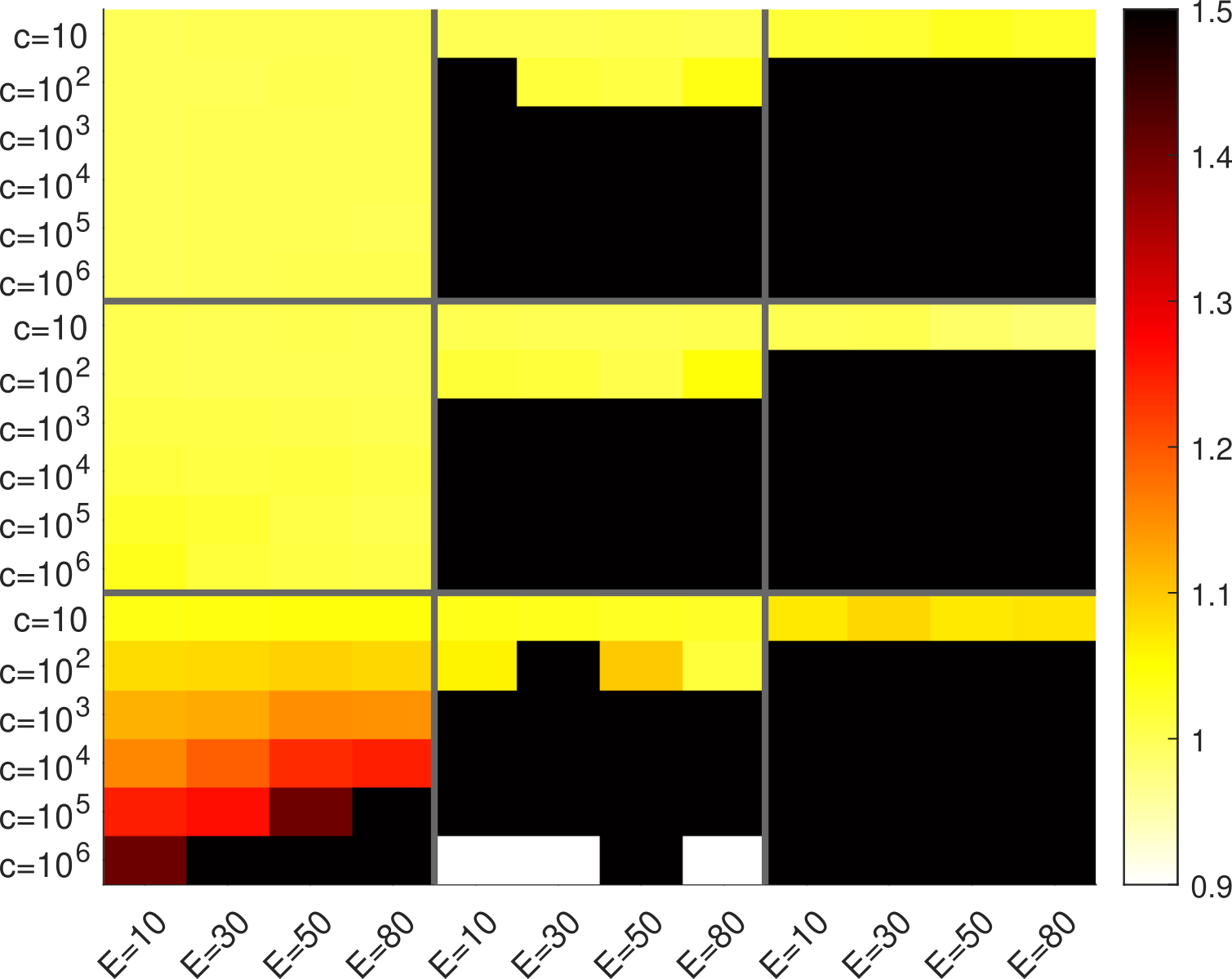}
    \includegraphics[width=0.48\columnwidth]{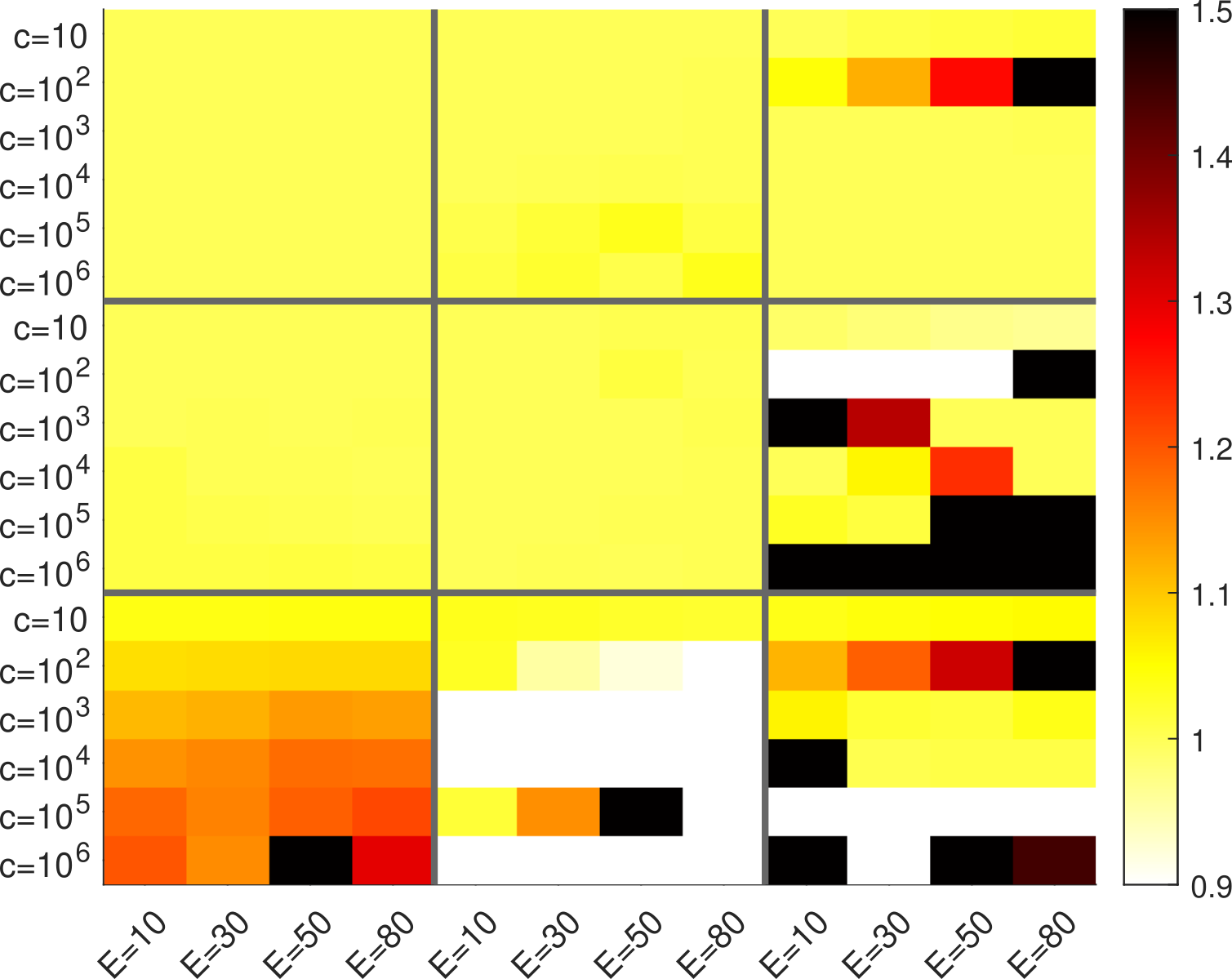}
    \caption{Performance heat-maps of the attained \textbf{best} objective function values by the MIESs when normalized with respect to CPLEX, after recording 10 runs over each of the 216 instances: [LEFT] \texttt{mies}, [RIGHT] \texttt{cma-IH}. Each of the $3 \times 3$ blocks, which are placed according to the problem-gallery of Figure \ref{fig:landscapes-gallery}, captures 24 instances as defined by the $(c,E)$ pairs. The colormap is scaled in $[0.9,1.5]$, where white reflects outperformance over CPLEX, yellow reflects similar performance as CPLEX, and the darkening red toward black reflects underperformance and possible divergence of the MIESs. 
    Importantly, there exist 3 problem instances on which CPLEX fails and terminates without a solution (white pixels at \textbf{(TC-7)} with $c=10^6$ and $E= \{\ 10,30,80 \}\ $). }
    \label{fig:cmaMap_best64}
\end{figure}
\subsubsection{High-Level Comparison: Best Attained Values}\label{sec:bestRes}
Figure \ref{fig:cmaMap_best64} presents a performance heat-map of the two MIESs based on the attained \textbf{best} objective function values when normalized with respect to CPLEX. It depicts the normalized best values, out of 10 recorded runs, over each of the 216 instances. It is organized as a $3 \times 3$ gallery, whose blocks are placed according to the problem-gallery of Figure \ref{fig:landscapes-gallery}, each capturing 24 instances as defined by the $(c,E)$ pairs. 
There exist 3 problem instances on which CPLEX fails and terminates without a solution (white pixels at \textbf{(TC-7)} with $c=10^6$ and $E= \{\ 10,30,80 \}\ $).

The overall performance of the \texttt{mies} is evidently weak -- it is outperformed by CPLEX in $\approx$75\% of the cases while obtaining similar results on the other cases. 
It does not locate solutions on the 3 (TC-7)-instances where CPLEX fails, but manages to outperform CPLEX on 2 specific instances of (TC-5) at $c=10$. \\
The \texttt{cma-IH}, however, performs very well on the majority of the problems.
It outperforms CPLEX on 28 problem instances, which include locating solutions to the 3 instances of (TC-7). Interestingly, this outperformance always involved the HadEllipse Hessian in either the objective or constraint functions. It performs similarly to CPLEX on 108 instances, and underperforms on 80 instances. \\
Altogether, it is evident that \textbf{CPLEX delivers inferior results on 28 out of the 216 instances} ($\approx$13\%).

Another high-level conclusion concerns the constructed MIQCQP test-suite -- it appears to possess richness of problem characteristics and difficulties, which are controlled by the combination of the Hessian functions in use, as well as the conditioning and the constraint level altogether. 
Among all problems, the MIESs' performance on (TC-6) consistently degraded as the conditioning and the constraint level increased. 
\subsubsection{Statistics of the MIESs' Runs}\label{sec:medRes}
\begin{figure}
    \centering
    \includegraphics[width=0.48\columnwidth]{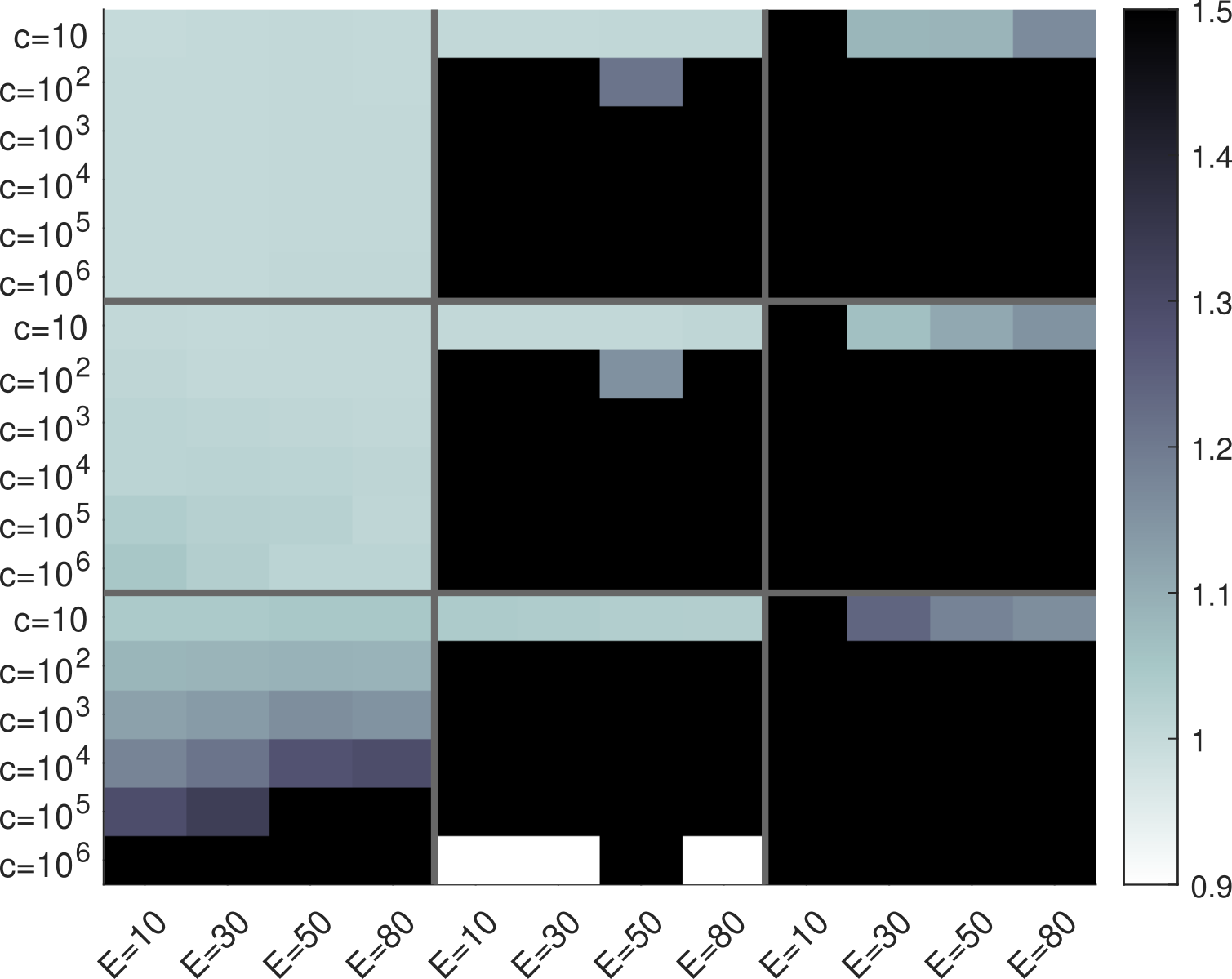}
    \includegraphics[width=0.48\columnwidth]{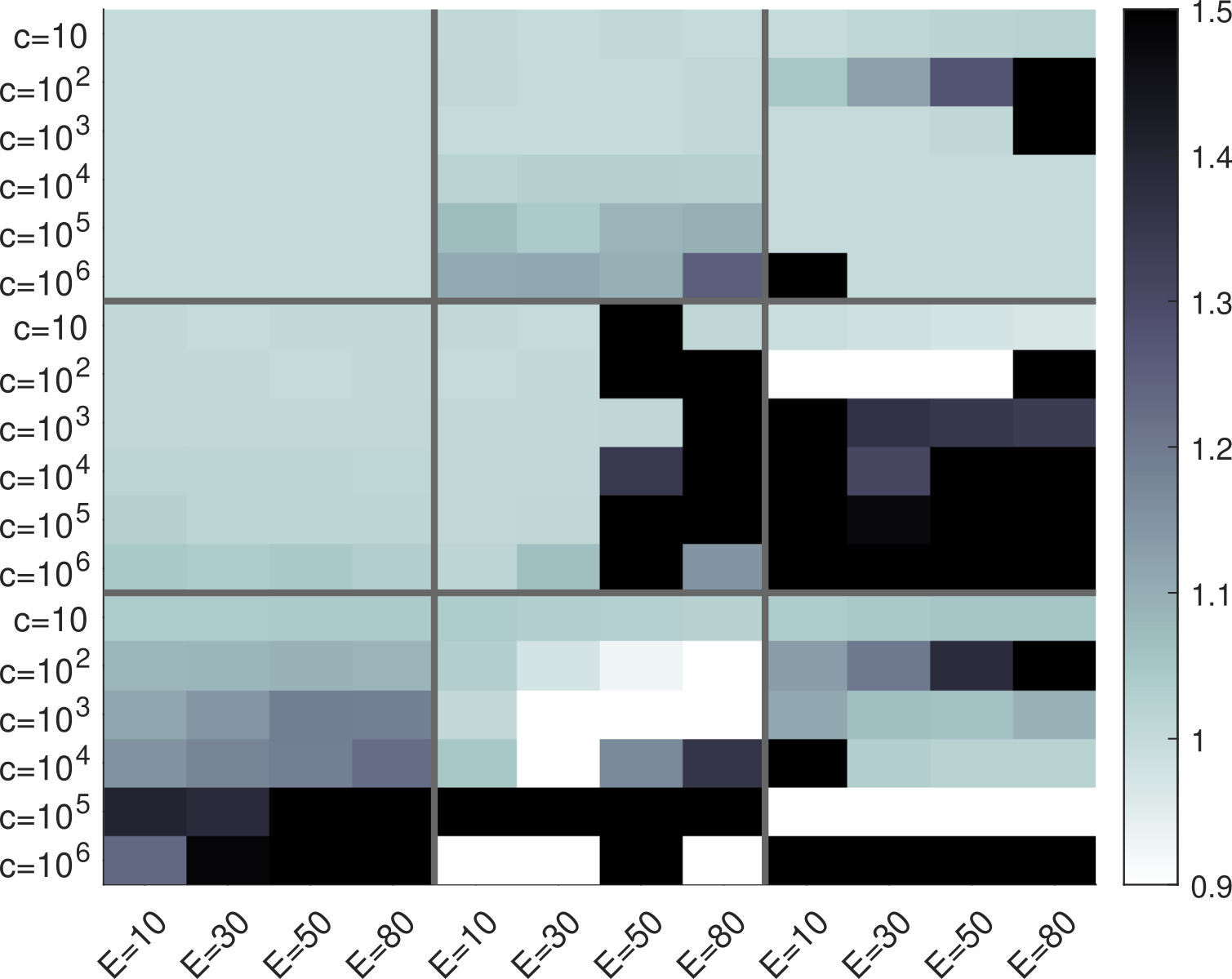}
    \caption{Performance heat-maps of the attained \textbf{median} objective function values by the MIESs when normalized with respect to CPLEX, recording performance of 10 runs over each of the 216 instances: [LEFT] \texttt{mies}, [RIGHT] \texttt{cma-IH}. The $3 \times 3$ gallery format is similar to Figure \ref{fig:cmaMap_best64}. The colormap is scaled in $[0.9,1.5]$, where white reflects outperformance over CPLEX, light grey/blue reflects similar performance as CPLEX, and black reflects underperformance and possible divergence of the MIESs.}
    \label{fig:cmaMap_all64}
\end{figure}

\begin{figure}
%\centering
\begin{tabular}{l | l }
    \hline 
    \scriptsize [A] \quad \textbf{(TC-0)} (Cigar+Cigar) & \scriptsize [B] \quad \textbf{(TC-6)} (HadEllipse+Cigar) \\ 
    \includegraphics[width=0.49\columnwidth,angle=0]{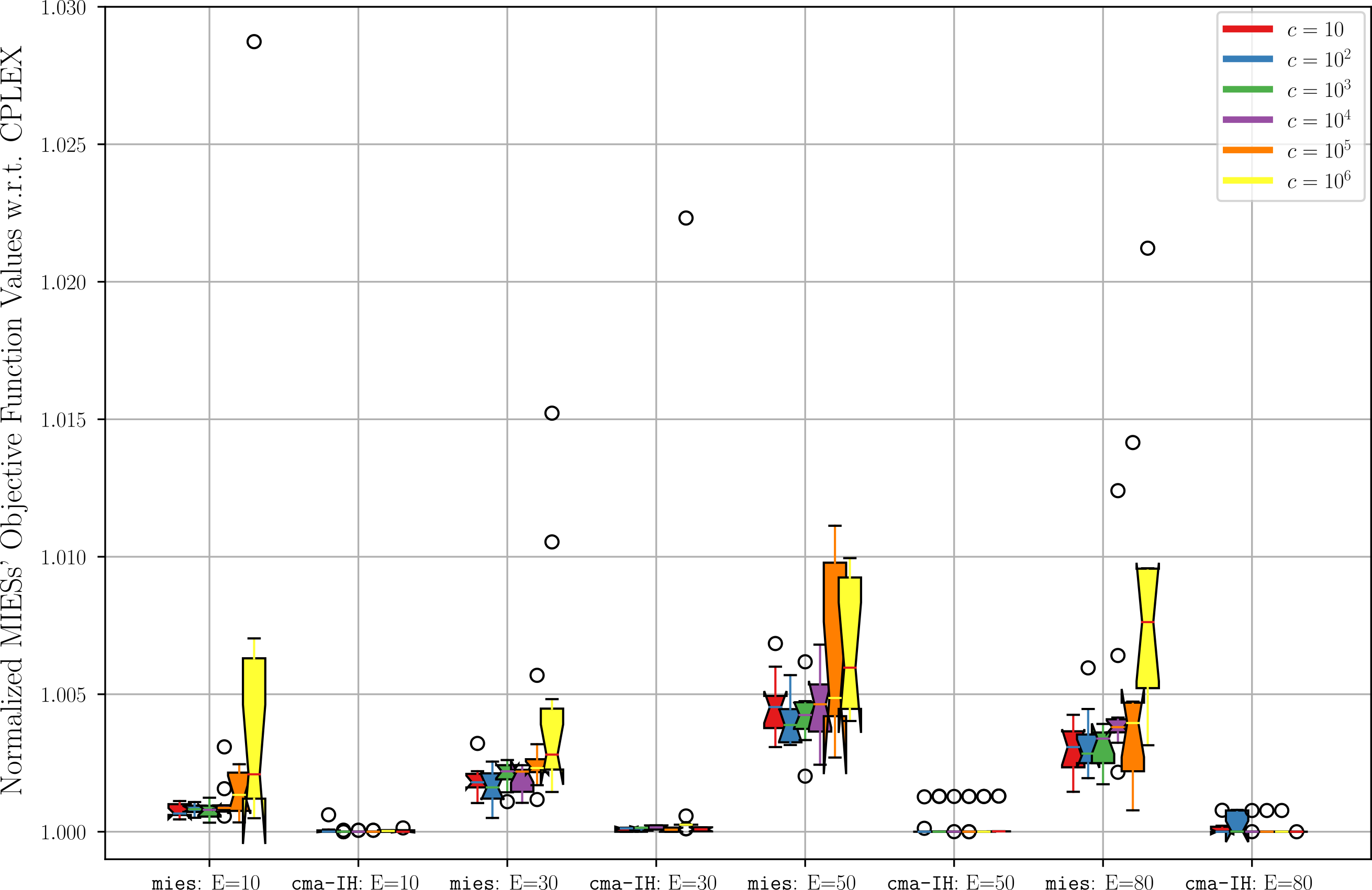} &
    \includegraphics[width=0.49\columnwidth,angle=0]{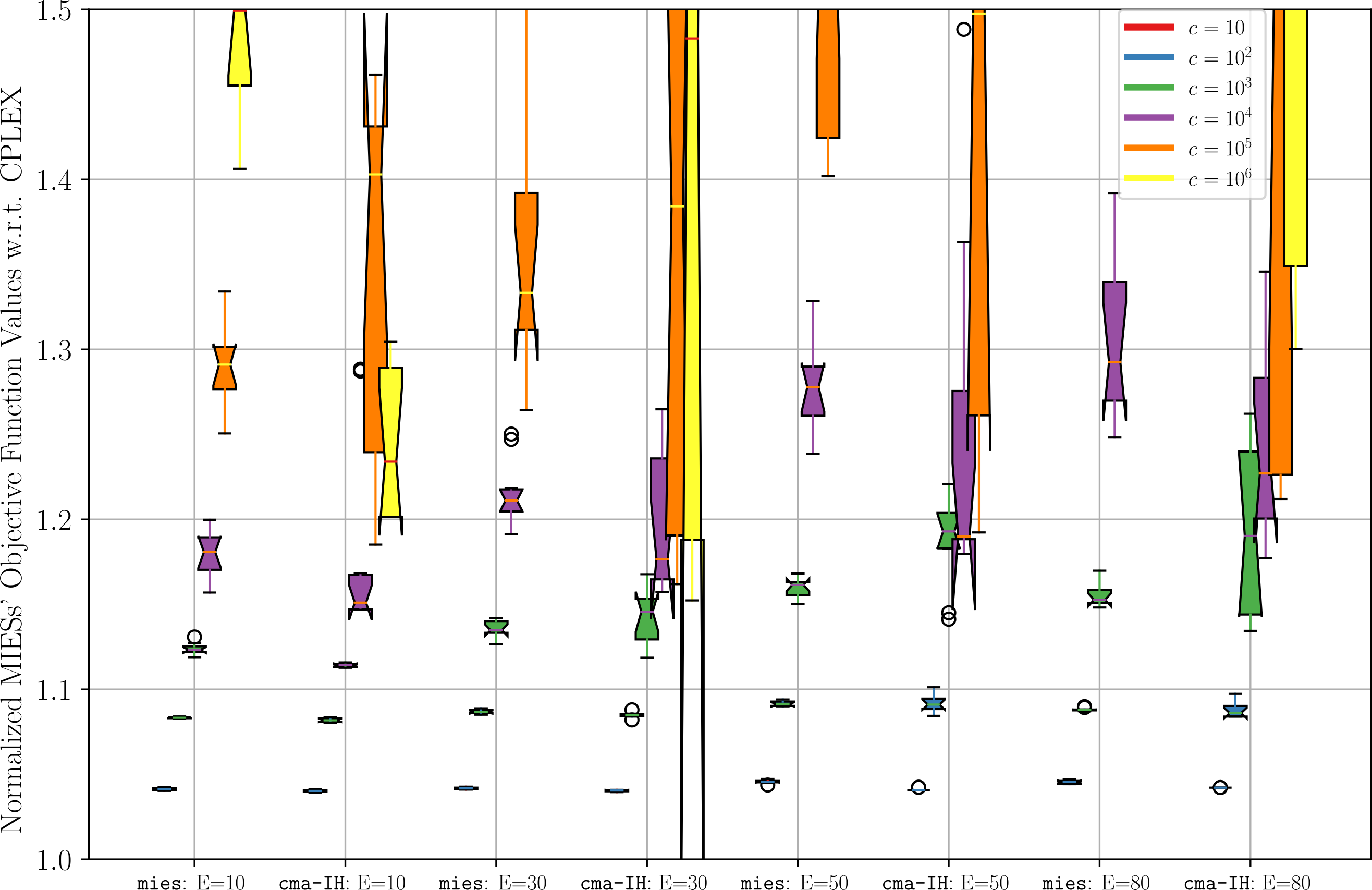} \\
    \hline
    \scriptsize [C] \quad \textbf{(TC-5)} (RotEllipse+HadEllipse) & \scriptsize [D] \quad \textbf{(TC-7)} (HadEllipse+RotEllipse) \\
    \includegraphics[width=0.49\columnwidth,angle=0]{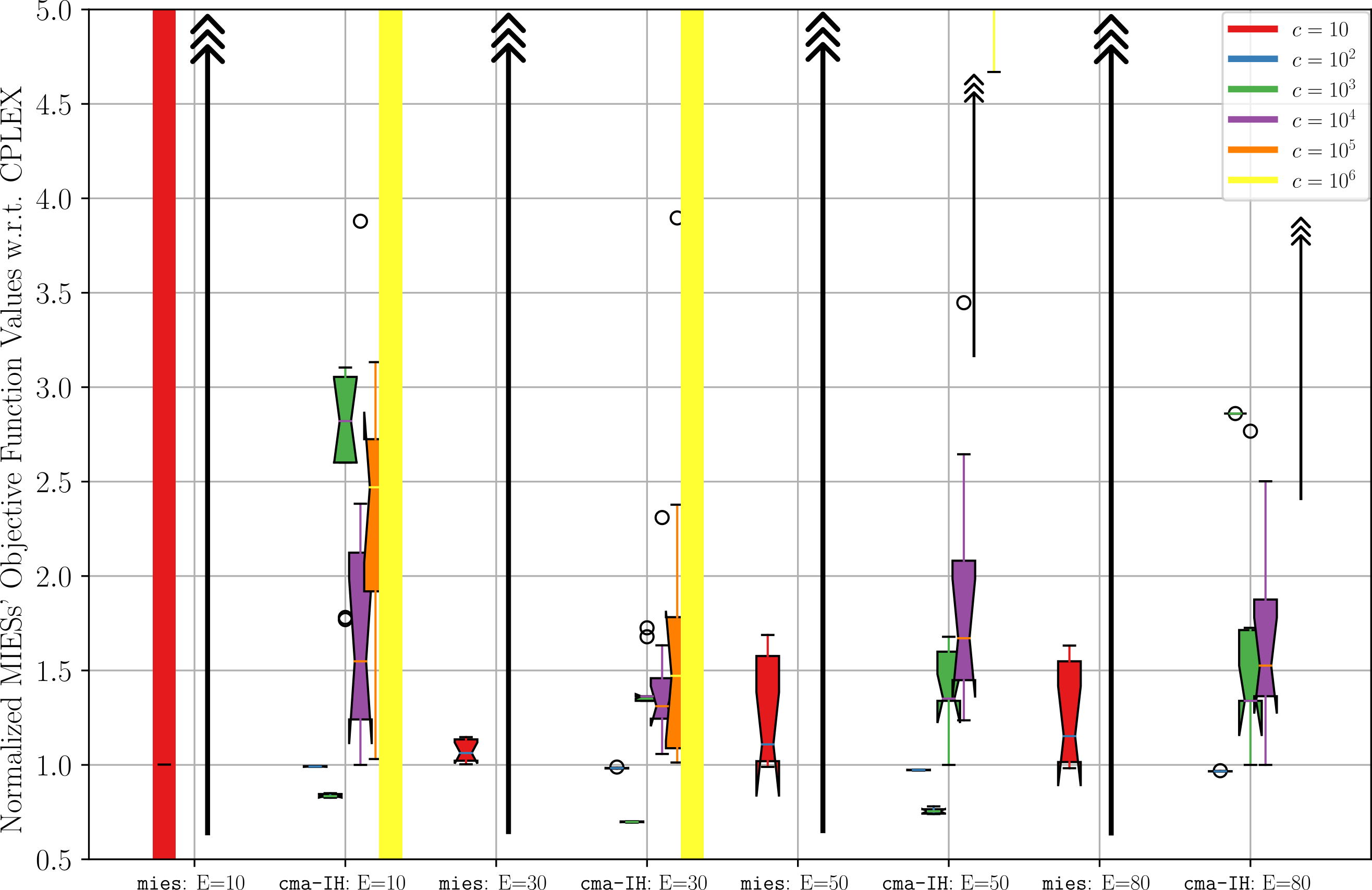} &
    \includegraphics[width=0.49\columnwidth,angle=0]{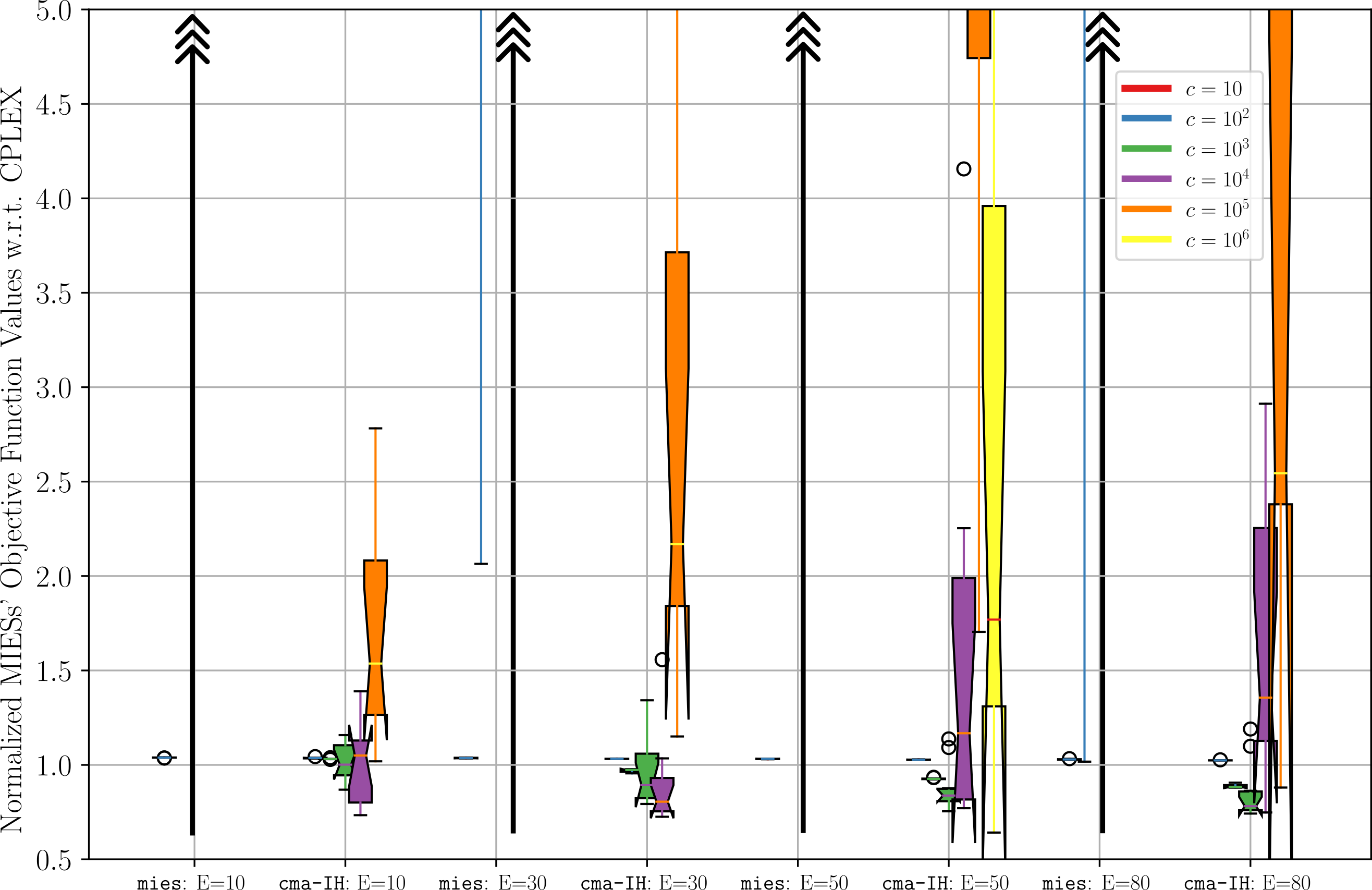} \\
    \hline
    \end{tabular}
    \caption{A gallery of statistical boxplots of the MIESs' objective function values normalized with respect to CPLEX per 4 problems.
Multi-arrow heads represent boxplots with extreme magnitudes (indicative of high degree of divergence). Mind should be given to the $y$-scale, which is not uniform across the gallery. Outliers are depicted as black circles. The boxplots are group-organized according to the constraint level $E$ (see $x$-axis' ticks) and group-colored according to the conditioning $c$. For each constraint level $E$, the \texttt{mies} groups are juxtaposed on the left-hand-side to the \texttt{cma-IH} groups.\label{fig:boxplotsObjF}}
\end{figure}

Figure \ref{fig:cmaMap_all64} presents a performance heat-map of the two MIESs based on the attained \textbf{median} objective function values when normalized with respect to CPLEX. It depicts the normalized median values of 10 recorded runs over each of the 216 instances. It is organized in the same $3 \times 3$ gallery format. 
A thorough analysis of each block was carried out in higher resolution using statistical boxplots, which are available online.\footnote{Raw datasets and boxplots will be provided upon request.}
As an example, we present 4 statistical boxplots for (TC-0) and (TC-3) in Figure \ref{fig:boxplotsObjF}, where the populations of the objective function values obtained by the stochastic MIESs are normalized with respect to the deterministic CPLEX.
%(that is, the value of 1.0 always represents the result attained by CPLEX): the separable Cigar as the constraint function [(TC-0),(TC-1)] (depicted in Fig.~\ref{fig:boxplotsObjF}[A,B,E,F])  and the non-separable RotEllipse [(TC-2),(TC-3)] (depicted in Fig.~\ref{fig:boxplotsObjF}[C,D]) serving in that role.
The boxplots are group-organized according to the constraint level $E$ (see $x$-axis' ticks) and group-colored according to the 6 conditioning levels (see legends), with 24 boxplots per ES. Next, we delve into the details.
%%%%%% =====> FIGURES TO UNCOMMENT
\normalsize 
\paragraph{\texttt{mies}}
The \texttt{mies} performed well when the Cigar function played the role of the constraint: smooth performance on (TC-0) (usually attaining objective function values within 1\% away from the global optima; see Figure \ref{fig:boxplotsObjF}[A]), good performance on (TC-3), and exhibited fine performance on (TC-6), which degraded as the conditioning increased (see Figure \ref{fig:boxplotsObjF}[B]). Otherwise, it performed poorly, or diverged, on the remaining 6 problem cases, except when the conditioning was set as low as $c=10$.

\paragraph{\texttt{cma-IH}}
The \texttt{cma-IH} usually performed very well on this test-suite also when considering its averaged-case behavior. %and outperformed the CPLEX over 20 problem instances . \\
This ES exhibited a similar pattern of performance as the \texttt{mies} when the Cigar function played the role of the constraint [left column of blocks: (TC-0), (TC-3), (TC-6)]. 
Notably, the \texttt{cma-IH} often terminated on (TC-0) within 1min, yielding altogether a population of 10 accurate runs much faster than the CPLEX, which always terminated due to \td{\texttt{time-out}}.\\ 
Considering the other problems, it usually performed better than the \texttt{mies}.
When the RotEllipse function served as the constraint function (middle column of blocks), its performance was usually good on low and medium conditioning, but the constraint level also played a role.  
On (TC-1), similar to (TC-6), the performance was excellent on the low conditioning instances $c=\{10,100\}$, and has consistently declined with the increasing conditioning to objective function gaps in the order of 10\% (when considering the median run; see the median bars of $c=\{10^5,10^6\}$ in Fig.~\ref{fig:boxplotsObjF}[B]). 
On (TC-4), however, the decline in performance was not due to the increasing conditioning, but rather due to the constraint level $E$. Interestingly, (TC-4) constituted the easiest problem for the CPLEX solver (instant computation with a status \an{\texttt{tolerance}}), except for the low conditioning of $c=10$ (where it terminates due to \td{\texttt{time-out}}).\\
The \texttt{cma-IH} performed extremely well on (TC-7), except for certain instances (mainly involving the conditioning $c=10^5$). It outperformed CPLEX on multiple instances. \\
Finally, the analysis of the right column of blocks [(TC-2), (TC-5), (TC-8)], which share the HadEllipse as their constraint function, unveils a high diversity of performance comparisons. 
(TC-2) appears to be an easy problem for the \texttt{cma-IH}, except for certain instances. 
(TC-5), on the other hand, seems like a hard challenge except for the low conditioning instances at $c=\{\ 10,100 \}\ $. Interestingly, the \texttt{cma-IH} outperforms CPLEX over the $c=100$ instance, but otherwise it is outperformed by it over all instances when $c>100$. Lastly, the (TC-8) problem holds all sorts of performance comparisons -- increasing difficulty per $c=100$, strong outperformance and underperformance at $c=10^5$ and $c=10^6$, respectively, and otherwise fine performance of the \texttt{cma-IH}.

\subsection{Unboundedness Amplification via Significant Translation of Optima}\label{sec:translation}
\begin{figure}
    \centering
    \includegraphics[width=0.48\columnwidth]{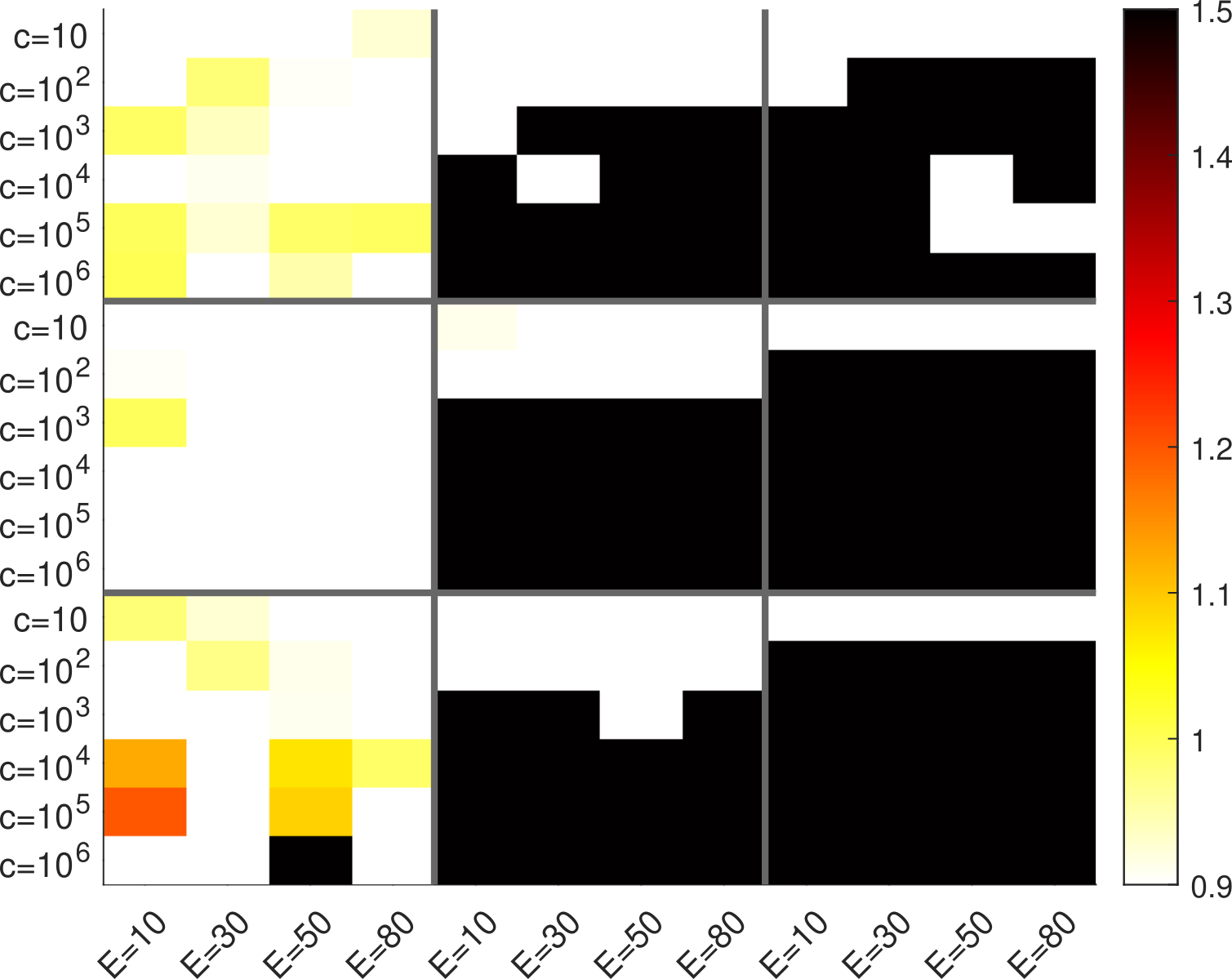}
    \includegraphics[width=0.48\columnwidth]{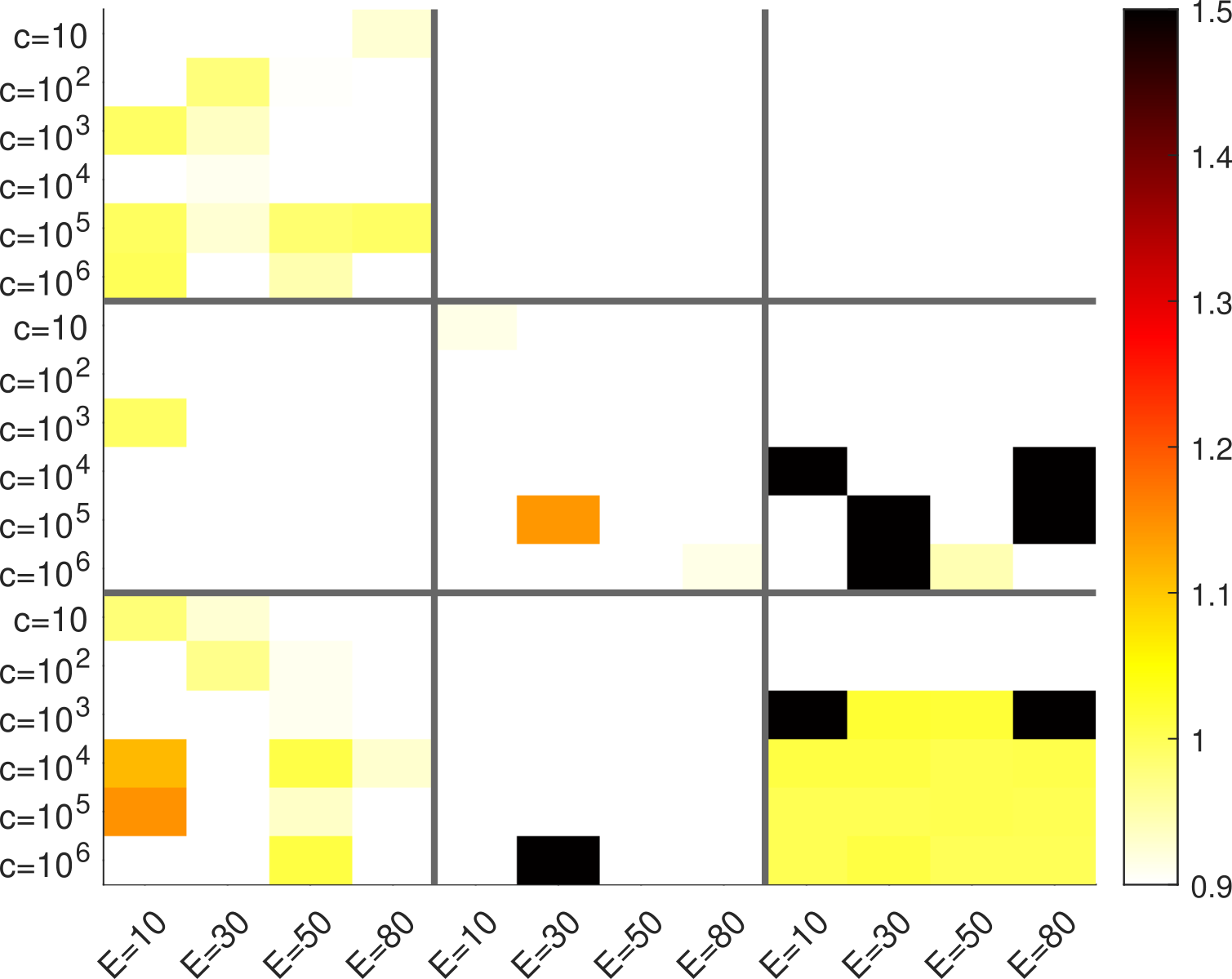}
    \caption{Heat-maps of the attained \textbf{best} objective function values by the MIESs on the \textbf{translated optima problem}, depicting normalized best attained value of 10 runs over each of the 216 instances (scaled in $[0.9,1.5]$; similar interpretation as in Figure \ref{fig:cmaMap_best64}): [LEFT] \texttt{mies}, [RIGHT] \texttt{cma-IH}.}
    \label{fig:Tmap}
\end{figure}

By further addressing the unboundedness research question, we are also interested in exploring the impact of translating the optima farther away from the origin.
Additional experiments were run to investigate the impact of a significant translation on the behavior of CPLEX and the MIESs. 
In practice, we account for the effect of shifting the vectors $\{ \vec{\xi}_0,\vec{\xi}_1 \}$ \eqref{eq:xi} as follows,
\[ 
\tilde{\vec{\xi}_i}:=\vec{\xi}_i + 10^4 \quad i=0,1,
\]
by conducting the same set of simulations prescribed in Section \ref{sec:setup}. \\
Similarly to the results reported in Section \ref{sec:benchmark}, numerical observations for the translated scenario are presented in Figure \ref{fig:Tmap}.
%and \ref{fig:Tcma}[LEFT] and .\\ %where notably, the MIESs occasionally 
\textbf{A groundbreaking finding is that the performance trend of \texttt{cma-IH} versus CPLEX is flipped }-- CPLEX becomes the inferior method on the translated use-case in 184 out of the 216 problem instances (underperformance of $\approx$85\%; see Figure \ref{fig:Tmap}[RIGHT]) unlike its status as the superior method when no translation occurred (underperformance in only $\approx$13\% of the cases; Figure \ref{fig:cmaMap_best64}[RIGHT]). 
The \texttt{mies} now outperforms CPLEX in $\approx$49\% of the cases.
Part of this behavior is due to the fact that CPLEX terminates more frequently without a solution. 
The \texttt{mies}, on the other hand, exhibits a consistent performance when solving the problems under translation -- its divergence on the challenging problem instances still holds (loosely speaking, the non-separable constraint functions' blocks remain ``black'', except for $c=10$), whereas it clearly outperforms CPLEX in its attainable problem instances (when the constraint function is the Cigar). 

\begin{figure}
    \centering
    \includegraphics[width=0.48\columnwidth]{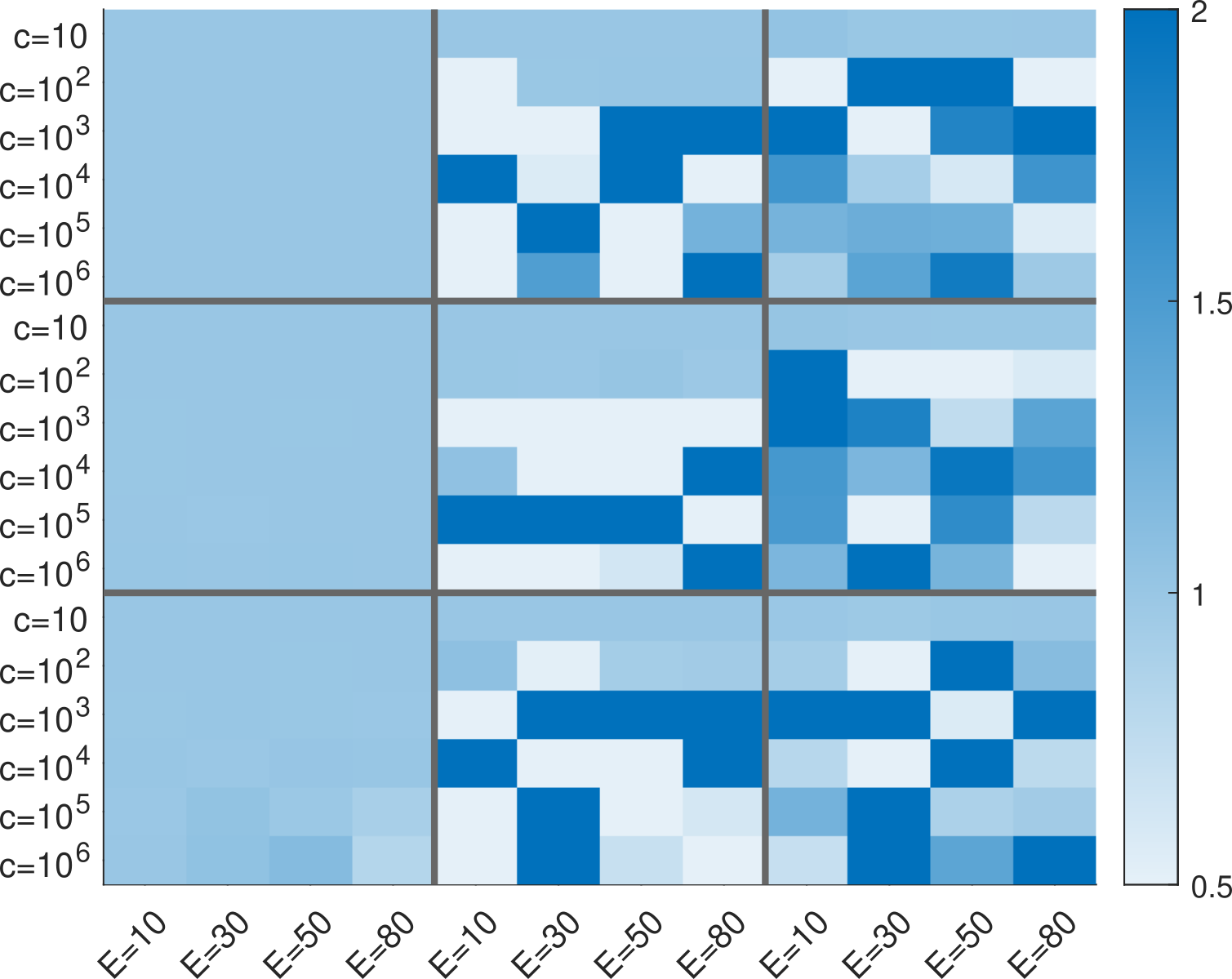}
    \includegraphics[width=0.48\columnwidth]{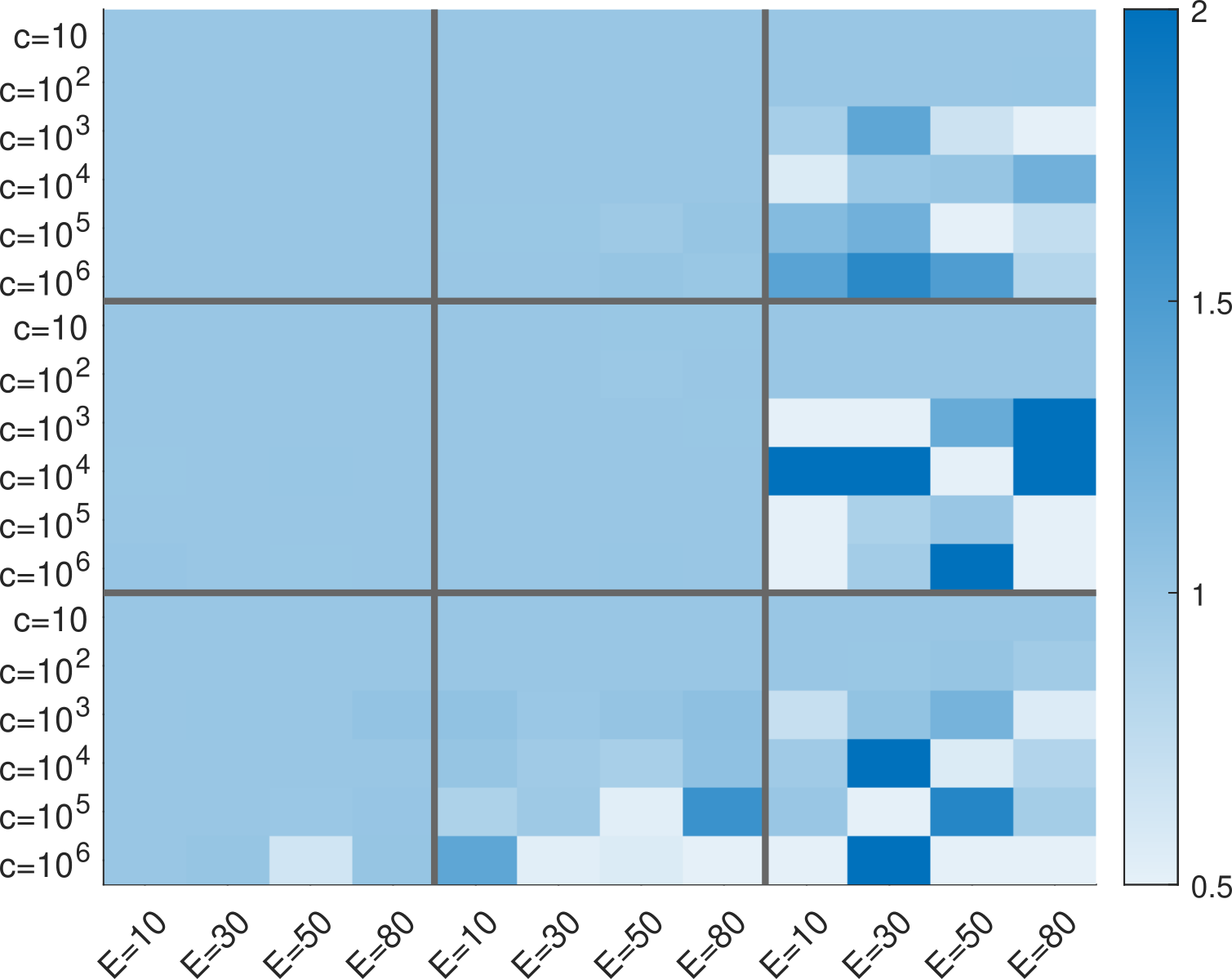}
    \caption{Heat-maps of the attained best objective function values by the MIESs on the \textbf{translated optima problem} when compared to the attained best values of the same strategy under no translation (scaled in $[0.5,2.0]$; light blue reflects similar performance, white reflects outperformance, and darker blue reflects underperformance): [LEFT] \texttt{mies}, [RIGHT] \texttt{cma-IH}.}
    \label{fig:Trelative}
\end{figure}

Furthermore, the MIESs' best runs under translation are compared to their own best runs under no translation on Figure \ref{fig:Trelative} (scaled in $[0.5,2.0]$; light blue reflects similar performance, white reflects outperformance of the translation case over the no-translation case, and darker blue reflects underperformance). 
It is evident that while most of the performance records remain similar upon introducing the translation | the \texttt{cma-IH} terminates in the exact same state in $\approx$65\% of the cases, whereas the remainder of the cases split between enhanced performance ($\approx$19\%) and degraded performance ($\approx$16\%). 
This behavior is to be expected when considering the self-adaptation principle \cite{Baeck-book,ESchapter2018}, which enables the ES to increase its step-sizes and to reach regions that are placed far-away from the starting points. Interestingly, the same principle may happen to assist in locating better optima in some of the cases, as if this long pathway with large step-sizes eventually enables a better convergence process, while hampering other cases.
At the same time, the \texttt{mies} maintains similar termination states in $\approx$41\% of the cases, whereas the remainder are qualitatively similar in ``divergence status'', as mentioned above.

\subsection{Integer Discrepancy}
\begin{figure}
    \centering
    \includegraphics[width=0.48\columnwidth]{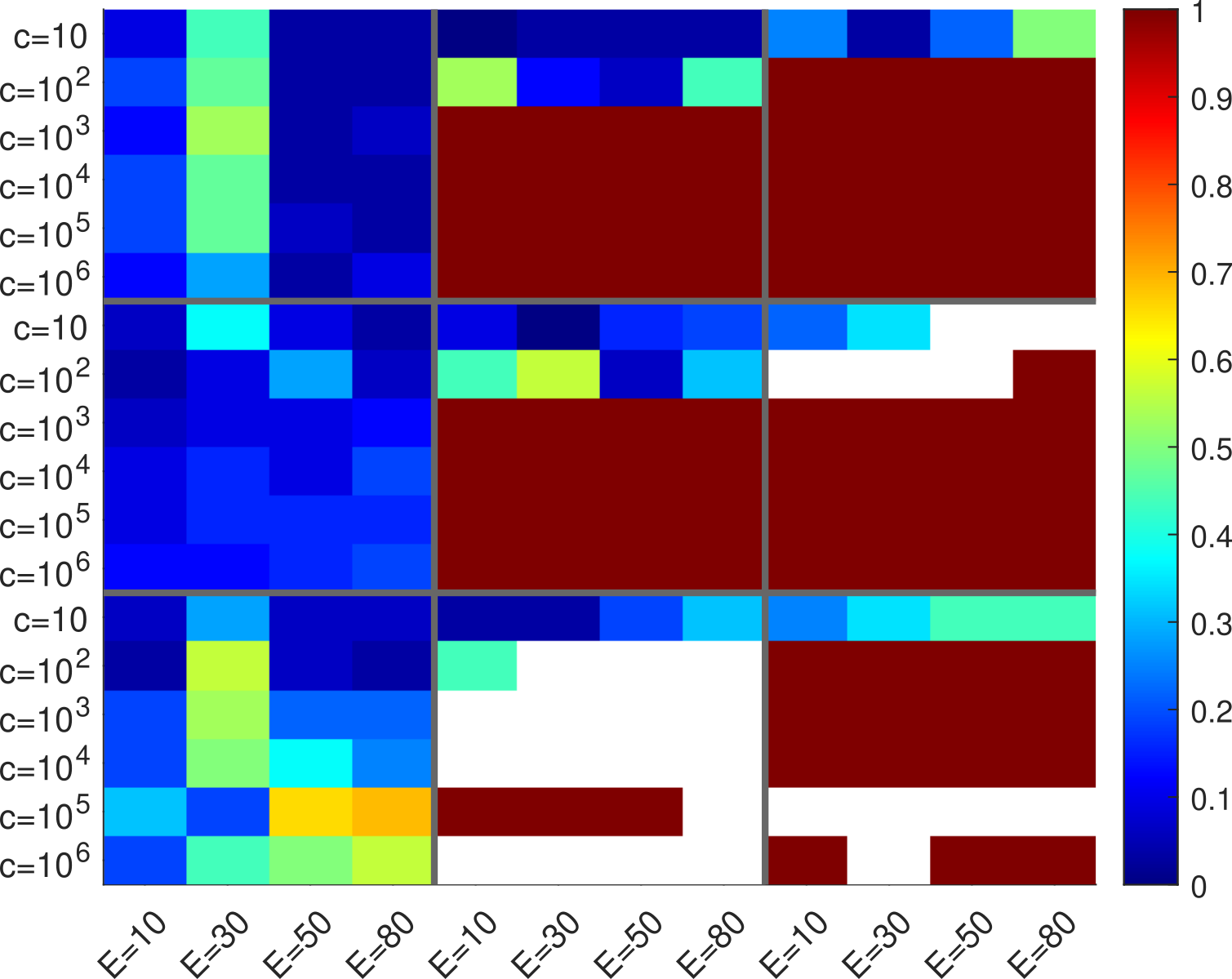}
    \includegraphics[width=0.48\columnwidth]{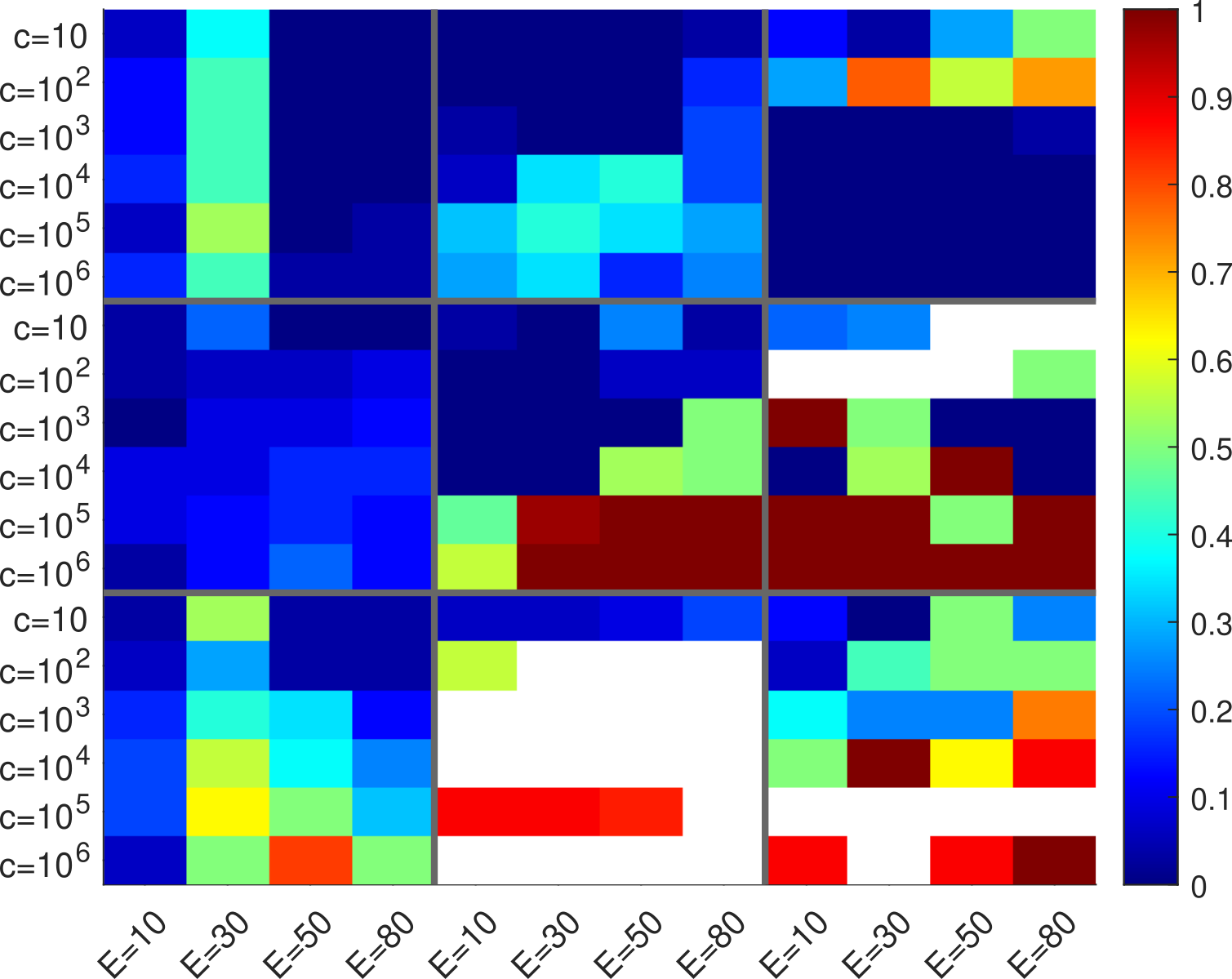}
    \caption{Heat-maps of the MIESs' normalized delta measure \eqref{eq:iER} when compared to the optimizers located by the CPLEX. Blue reflects similar integer vectors and red reflects total discrepancy: [LEFT] \texttt{mies}, [RIGHT] \texttt{cma-IH}.}
    \label{fig:DeltaIntMaps}
\end{figure}

When analyzing the numerical results in light of a comparison between the two MIESs, most of the times there is an advantage of the \texttt{cma-IH} with respect to the attained objective function values. The gap is occasionally negligible (e.g., below 1\% per (TC-0)), but sometimes prominent (either with statistical significance or without). 
%And yet, due to the mixed-integer nature, we are ought to assess the attained degree of precision in terms of the ``correct'' integers within the converged solutions. \\
And yet, due to the different mutation mechanisms of the two strategies, we would like to analyze the nature of the attained solutions in more detail and seek patterns in their behavior. \\
In order to numerically assess the precision of a candidate solution $\vec{\zeta}$ with respect to the integer optimizer $\vec{z}^{*}$ (i.e., the subset of integer decision variables at the global optimum), we introduce a normalized delta measure $\Delta_z \in [0,1]$ which counts the miscalculated variables belonging to the index set I:
\begin{align}\label{eq:iER}
\displaystyle \Delta_z\left(\vec{\zeta},\vec{z}^{*} \right) := & \displaystyle \frac {\# \left\{ i : \zeta_i \neq z_i^{*} ,~\forall i\in I \right\}}{n_z}~.
\end{align}
Care should be given when interpreting this delta measure on \textit{multimodal} landscapes that possess sub-optimal solutions of high-quality.\\
We will assume in what follows that CPLEX locates the global optimizer whenever it outperforms the MIESs.\\ 
Figure \ref{fig:DeltaIntMaps} presents heat-maps of the above delta measure when considering the best decision vectors attained by the MIESs when compared to the CPLEX vectors. The blank/white pixels correspond to the failures of the CPLEX.\\ 
The heat-maps unveil complex patterns, even when the MIESs' performance is solid as in (TC-0). 
We would like to delve deeper into the numerical data, and choose to restrict our investigation to the four cases involving only the Cigar and RotEllipse functions (also holding more confidence regarding the optimality of the CPLEX solutions). 

\begin{figure}
\caption{A gallery of statistical boxplots of the normalized delta measure $\Delta_z$ \eqref{eq:iER} uniformly scaled by $[0,1]$.
The boxplots are group-organized according to the constraint level $E$ and group-colored according to the conditioning $c$; outliers are depicted as black circles.
\label{fig:miER64}}
\centering
\begin{tabular}{c | c}
    \hline 
    \scriptsize [A] \quad \textbf{(TC-0)} (Cigar+Cigar) & \scriptsize [B] \quad \textbf{(TC-1)} (Cigar+RotEllipse) \\ 
    \includegraphics[width=0.49\columnwidth,angle=0]{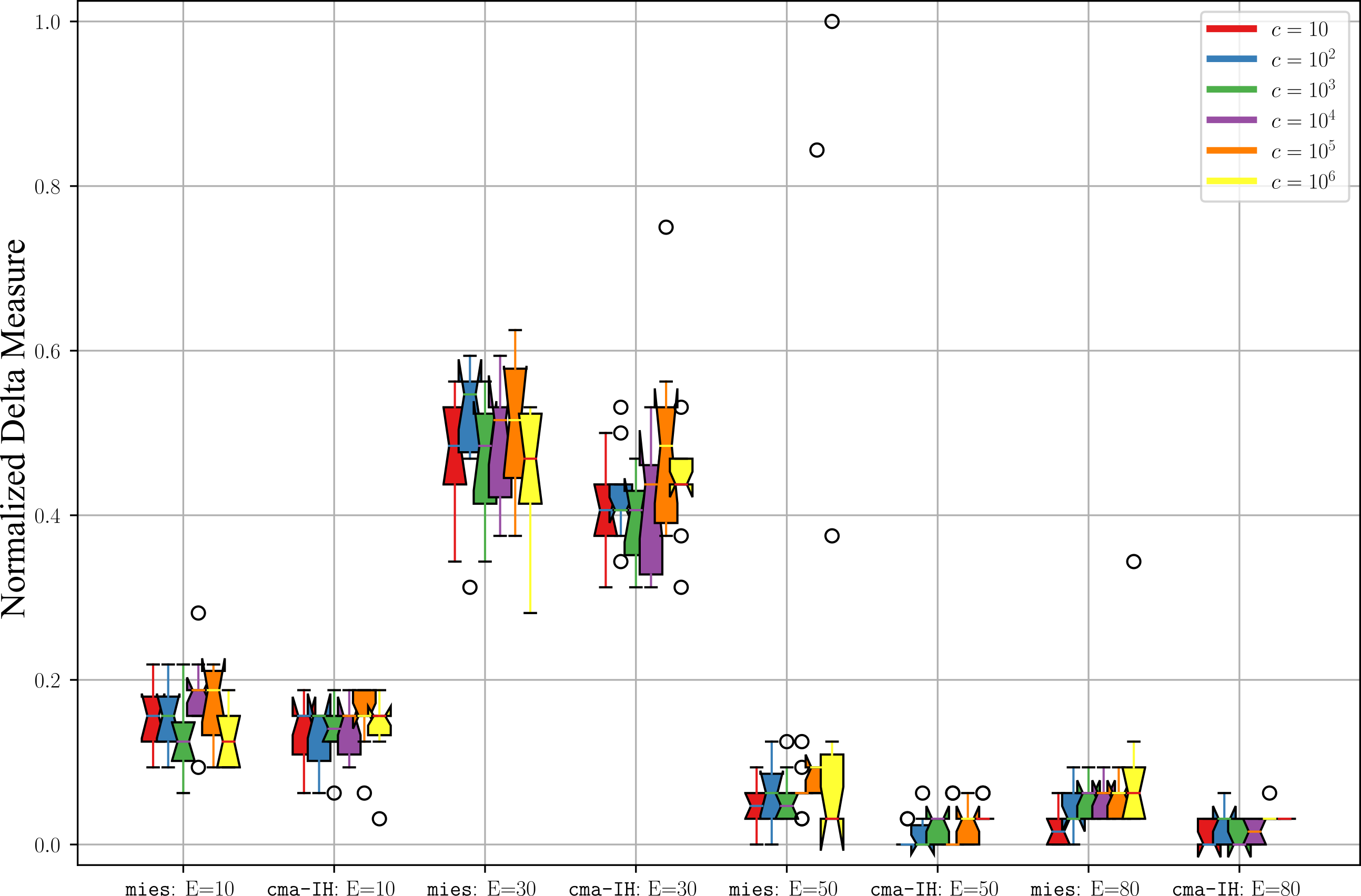} &
    \includegraphics[width=0.49\columnwidth,angle=0]{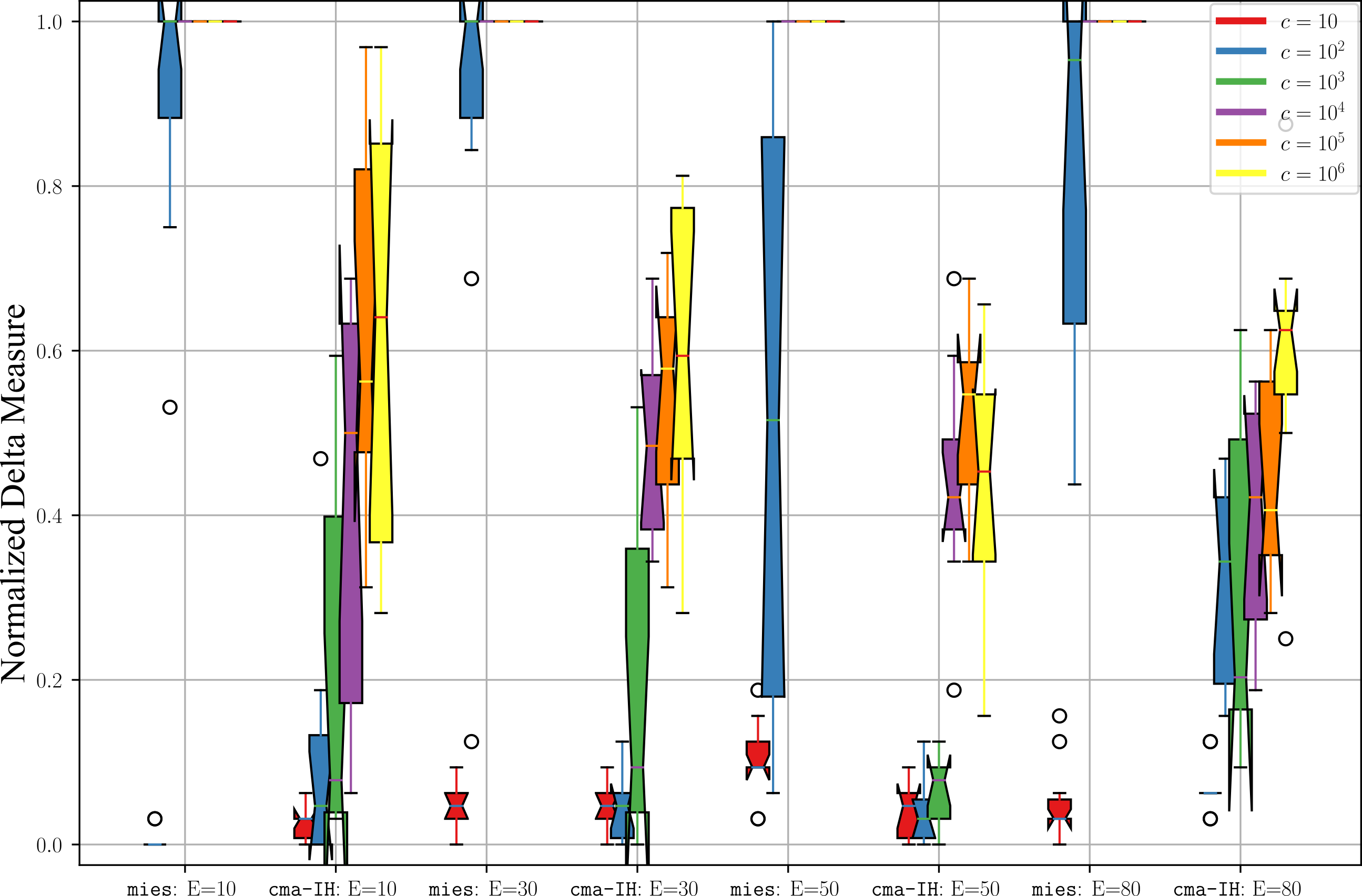} \\
    \hline
    \scriptsize [C] \quad \textbf{(TC-3)} (RotEllipse+Cigar) & \scriptsize [D] \quad \textbf{(TC-4)} (RotEllipse+RotEllipse) \\
     \includegraphics[width=0.49\columnwidth,angle=0]{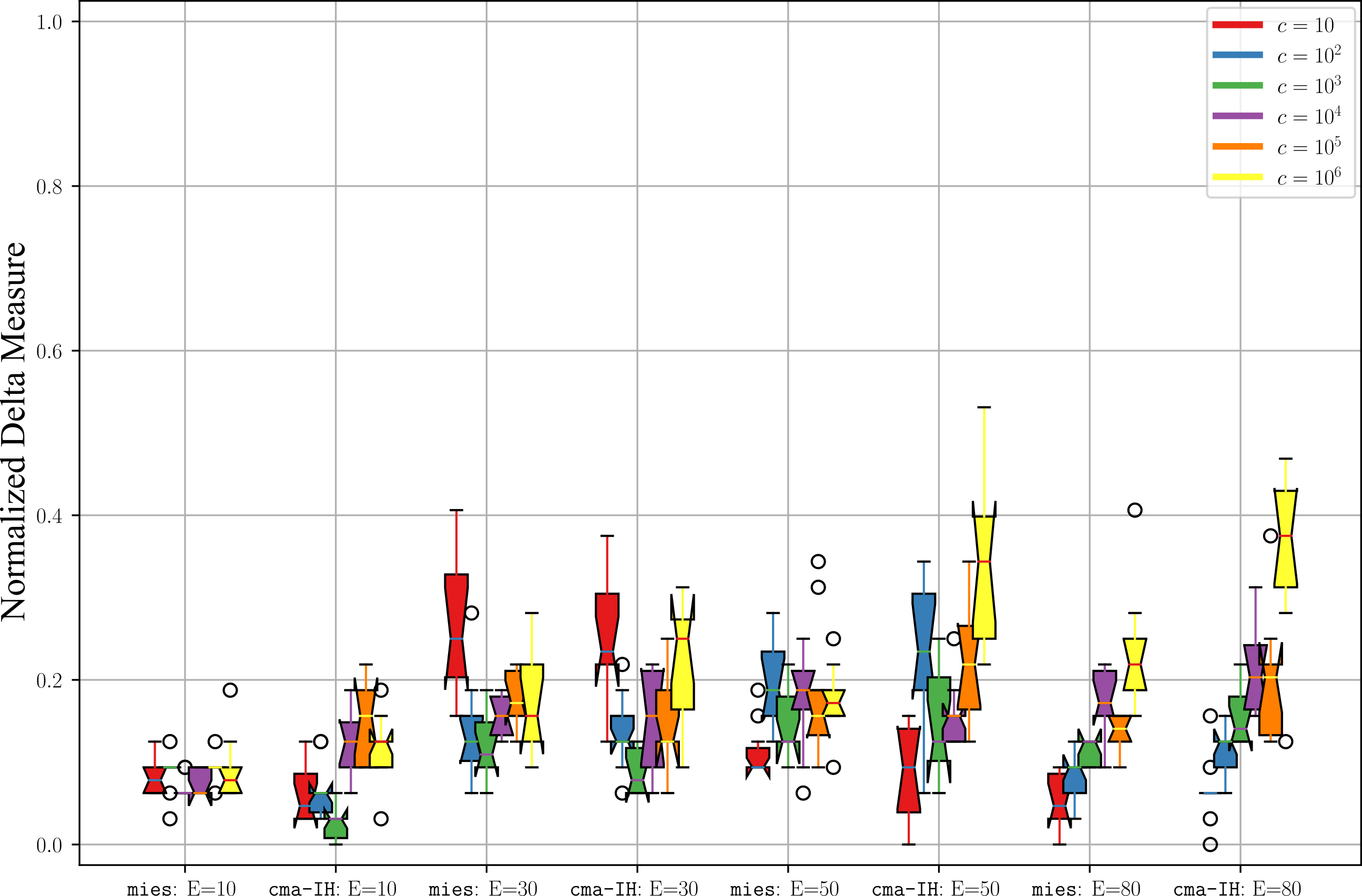} &
    \includegraphics[width=0.49\columnwidth,angle=0]{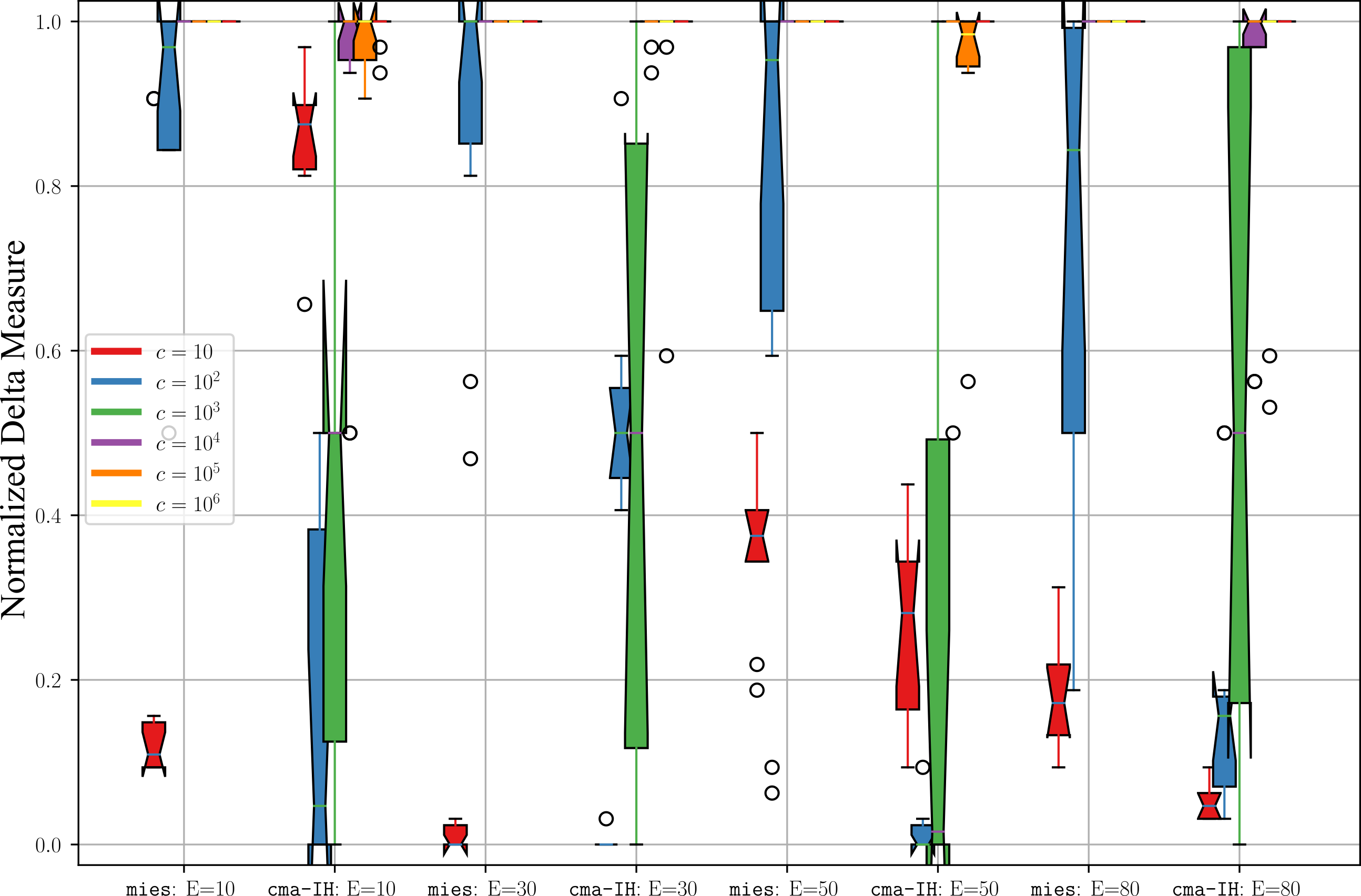} \\
    \hline 
    \end{tabular}
\end{figure}

Figure \ref{fig:miER64} presents the boxplots of the delta measure $\Delta_z$ for these four problems. 
%\td{TODO: elaborate on the error-rate - and add information on extreme cases of low/high errors}
As evident in Figure\ref{fig:miER64}[A], the two MIESs obtain similar delta measures for (TC-0). 
The constraint level of $E=30$ seems more challenging.
Per (TC-1), the delta measure pattern is equivalent to the objective function values' pattern, that is, the \texttt{mies} obtained low delta values only on $c=10$, while the \texttt{cma-IH} usually obtained low values that deteriorated (increased) as the conditioning increased (see Figure\ref{fig:miER64}[B]).
Next, when examining (TC-3), this trend of similar delta values is kept until the high conditioning of $c=\{10^5,10^6\}$, on which the trend flips and the \texttt{mies} obtains lower delta values with statistical significance (see Figure\ref{fig:miER64}[C]). 
Finally, (TC-4) exhibits a complex pattern of this delta measure. 
The \texttt{cma-IH} accomplishes fine objective function values (see Figure\ref{fig:boxplotsObjF}[D]), and yet its located solutions reside far away from the solutions located by the CPLEX (see Figure\ref{fig:miER64}[D]). Interestingly, solutions obtained by the \texttt{mies} are located within moderate delta values from those optima whenever it is successful ($c=10$). 
We hypothesize that this pattern is rooted in the \textit{multimodality} of the integer sub-space of this search landscape. We will address it in Section \ref{sec:RRmultimodality}.
We conclude that the advantage in the objective function attainment of the \texttt{cma-IH} over the \texttt{mies}, when they both converge, is rooted in its ability to accomplish better precision for the continuous part of the optimizers. The \texttt{mies} seems to identify the correct integers of the optimizers for (TC-0) and (TC-3), but to lack precision of the continuous part on these two test-cases. 
Altogether, the two MIESs seem competitive in addressing the unbounded integer search, but the \texttt{mies} seems to accomplish higher precision whenever it is able to tackle the problem. 
We will discuss it further in Section \ref{sec:summary}.

%%%============> FIGURES TO COMMENT
\normalsize

\subsection{All-Integer: Additional IQCQP Benchmarking}\label{sec:iqcqp}
\begin{figure}
\caption{Additional problem-solving of IQCQP (all-integer) at $D=64$, depicting the outcome only for (TC-1) (Cigar+RotEllipse): [LEFT] statistical boxplots of the MIESs' attained objective function values when normalized with respect to CPLEX, [RIGHT] the delta measure $\Delta_z$ \eqref{eq:iER} per the obtained populations of solutions (decision vectors). 
Multi-arrow heads represent boxplots with extreme magnitudes. 
Reference should be made to the mixed-integer problem-solving shown in Figure\ref{fig:cmaMap_all64} and Figure\ref{fig:miER64}[B].\label{fig:iqcqp64}}
    \centering
    \includegraphics[width=0.49\columnwidth]{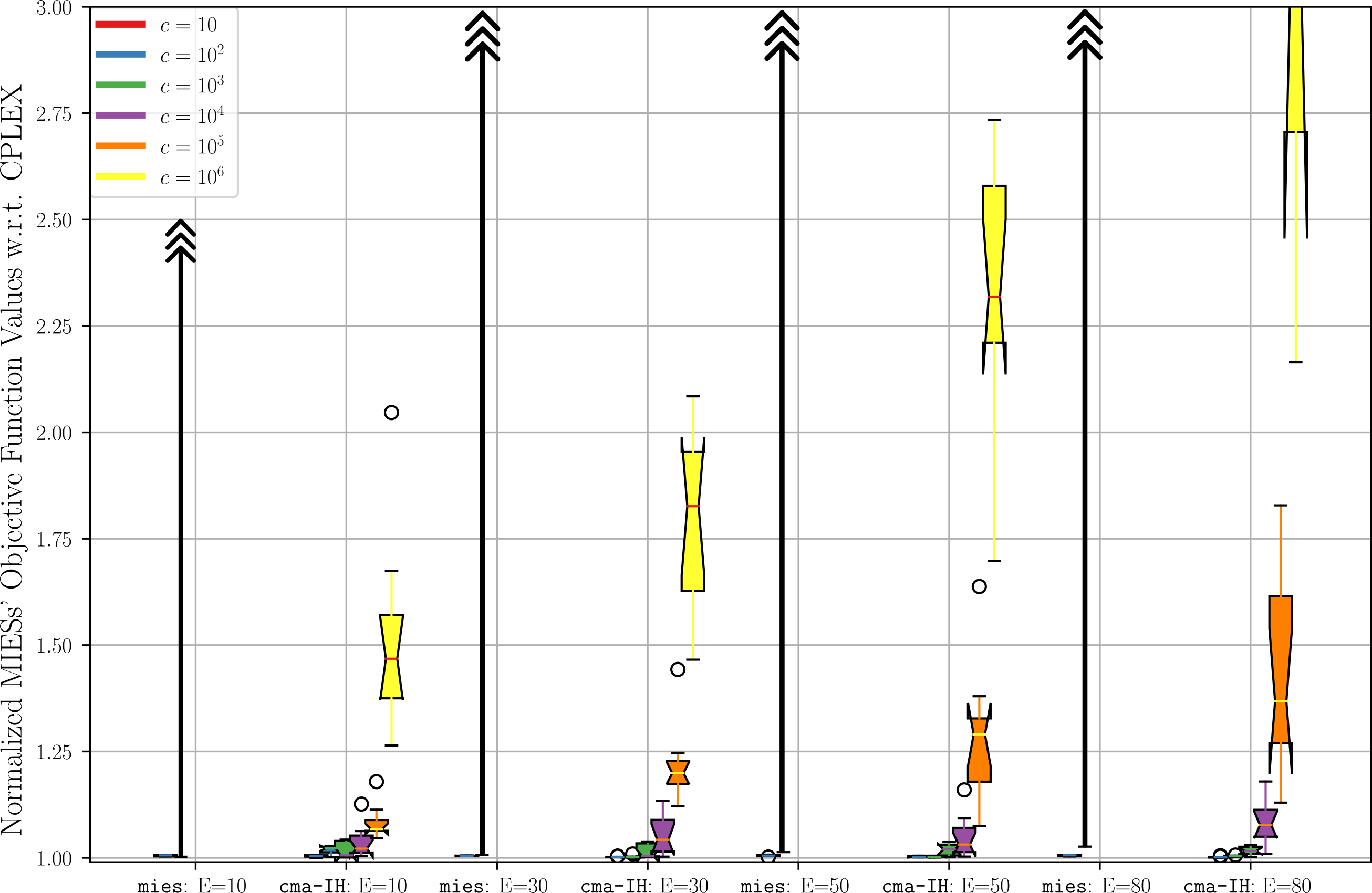} 
    \includegraphics[width=0.49\columnwidth]{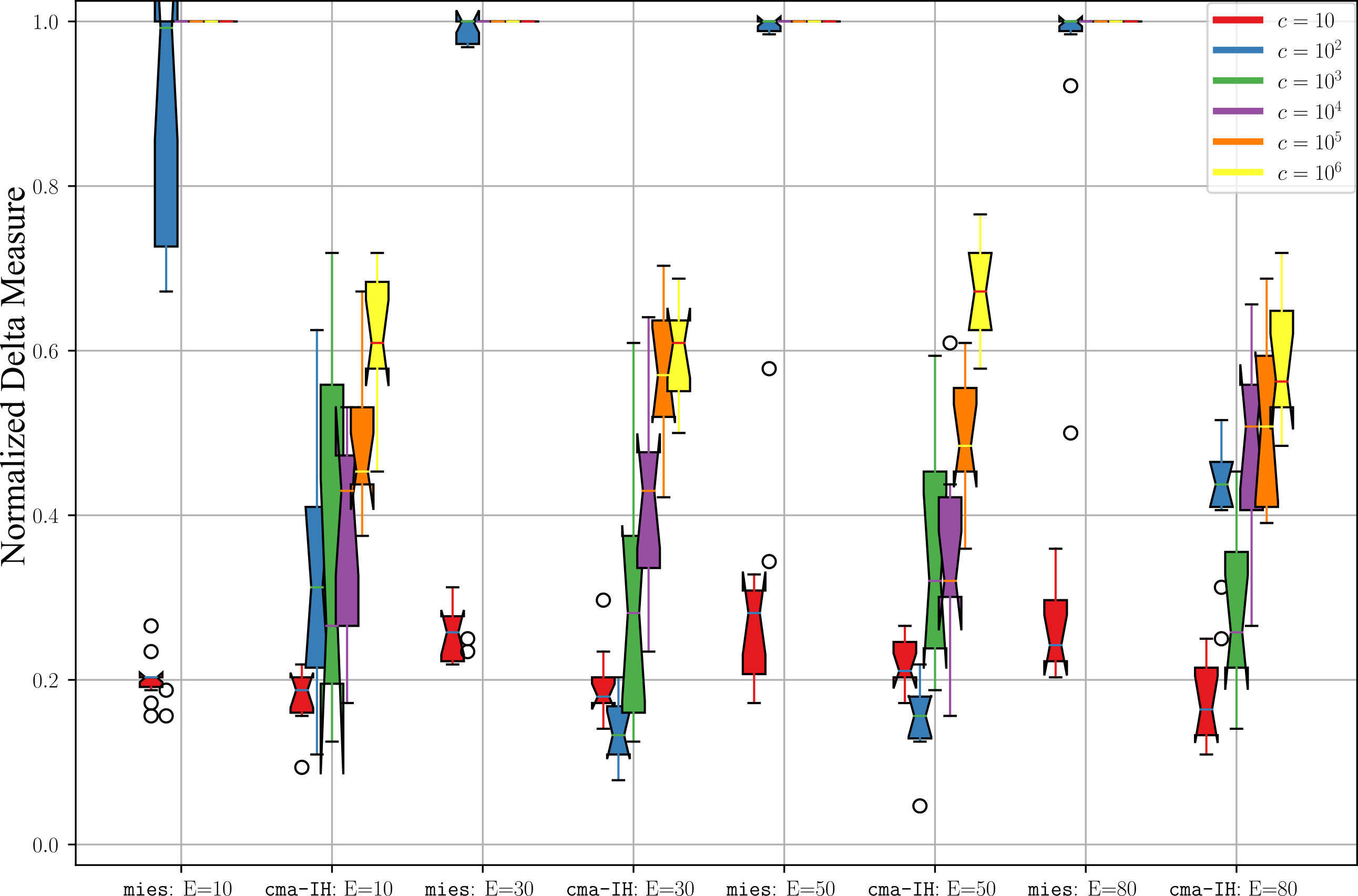}
\end{figure}

When turning to IQCQP instances (all-integer problems with equivalent specifications: $D=n_z=64,~n_r=0$), the observed statistical trends are consistent with the mixed-integer patterns shown earlier. It is evident that the two MIESs experience the same strengths and weaknesses when solving the all-integer problems. 
Due to space limitations, Figure \ref{fig:iqcqp64} presents the outcome per (TC-1) alone, comprising boxplots of the objective function values and the integer delta measure. Consistent with the mixed-integer trends shown in Figure\ref{fig:cmaMap_all64} and Figure\ref{fig:miER64}[B], the observed statistical patterns further reflect that the \texttt{cma-IH} is being challenged by this problem (mind should be given to the $y$-scale).
\subsection{Assessing Multimodality of the Integer RotEllipse Sub-Space}\label{sec:RRmultimodality}
\begin{figure}
    \centering
    \includegraphics[width=0.49\columnwidth,angle=0]{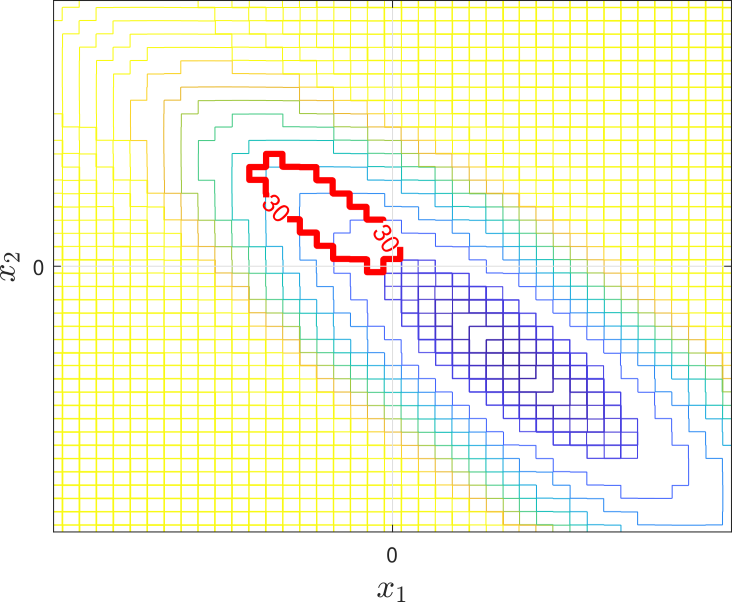} \includegraphics[width=0.49\columnwidth,angle=0]{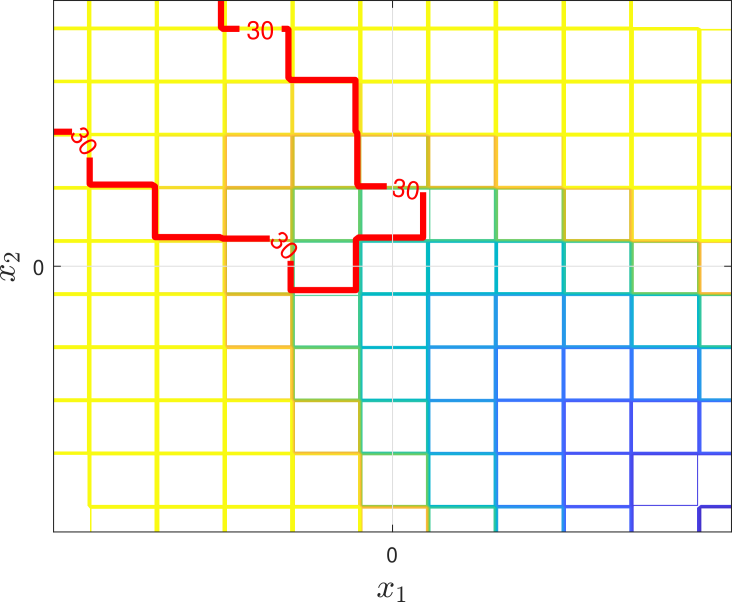} 
    \caption{Integer contours of the 2D RotEllipse+RotEllipse test-case (TC-4) at $c=10$ in two scales: $[-20,20]^2$ [LEFT] and  $[-5,5]^2$ [RIGHT]. The constraint function, a \textit{discrete feasibility ellipsoid} centered about (-4,4), is depicted by the thick red curve, and is set here to $E=30$. The objective function's ellipsoid is centered about (7,-7).} \label{fig:iRR10}
\end{figure}
\begin{figure}
    \centering
    \includegraphics[width=0.45\columnwidth,angle=0]{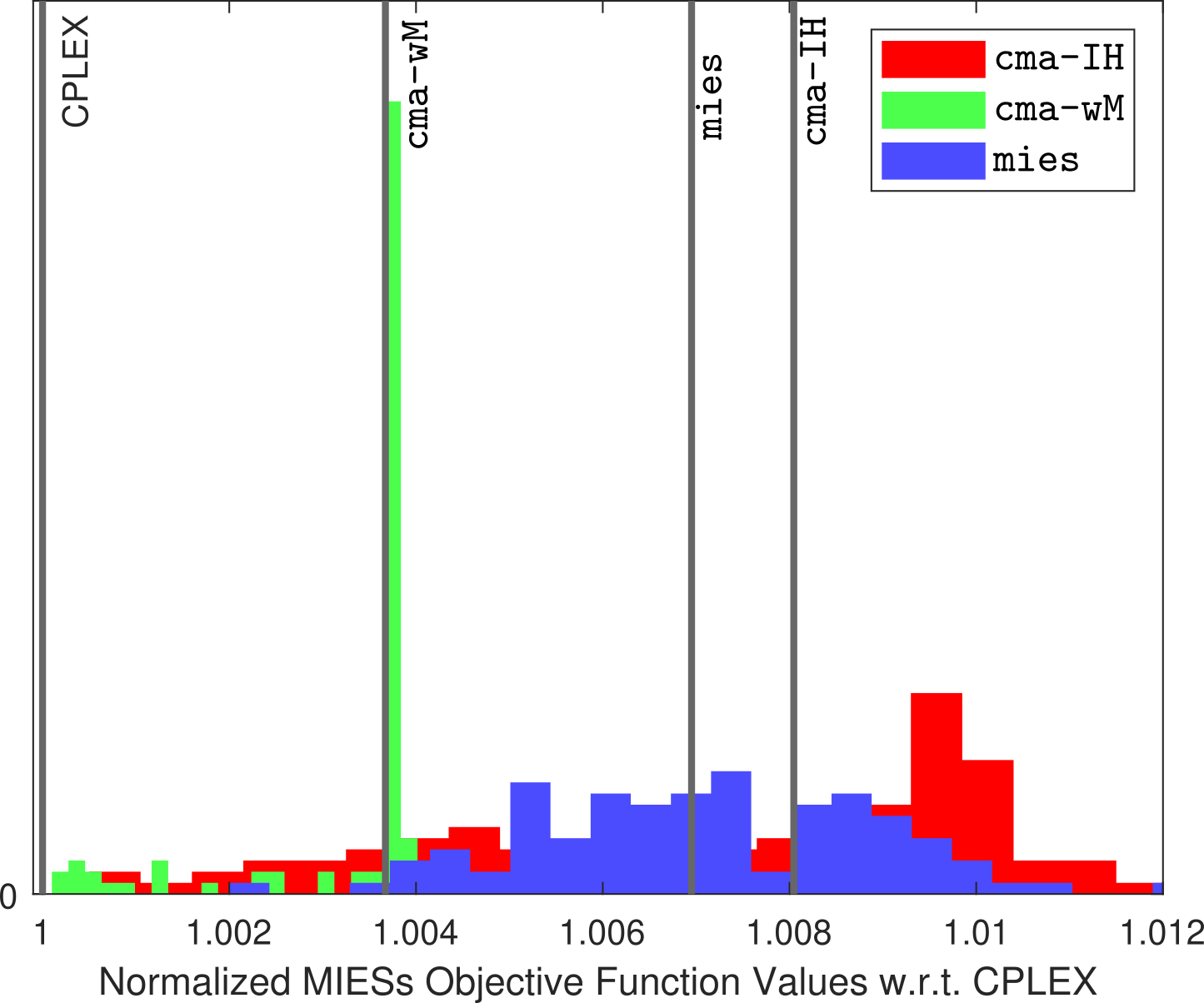} \includegraphics[width=0.52\columnwidth,angle=0]{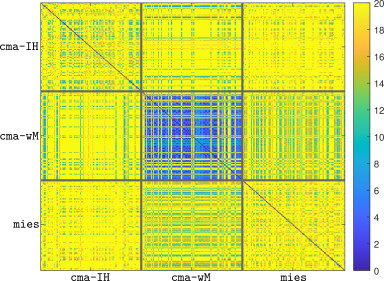} 
    \caption{Problem solving of the (bounded) all-integer 64D RotEllipse+RotEllipse test-case (TC-4) at $c=10$ constrained to $E=10$:  
    [LEFT] histogram of attained objective function values normalized with respect to CPLEX (represented by a vertical black line at $1.0$) by three MIESs (100 runs each): \texttt{cma-IH} (red), \texttt{cma-wM} (green), and \texttt{mies} (blue). The respective \textit{median} values are depicted as vertical black lines with adjacent labels. [RIGHT] Mutual $L_1$ distances among the attained solution vectors (colormap: $[0,20]$).} \label{fig:resIRR10}
\end{figure}

Following the observation that large gaps exist among the integer decision variables of the top solutions per (TC-4), we would like to address the hypothesis that its integer sub-space is multimodal. 
We begin by examining the plot of the reduced 2-dimensional test-case, depicted as Figure \ref{fig:iRR10}. 
Indeed, the intersection of the \textit{discrete feasibility ellipsoid} (in thick red) with the discrete level sets of the objective function generate altogether multipility of high-quality solutions. Importantly, this 2-dimensional test-case does not necessarily well represent the structure in high-dimensions, but it serves to equip us with intuition.\\
Next, we conducted extended simulations for the all-integer 64-dimensional RotEllipse+RotEllipse test-case for four concrete instances: $c=\{10,1000\}$ and $E=\{10,30\}$. Since our target here is to analyze the landscape, we focused on obtaining as many high-quality solutions as possible, regardless of the unboundedness question. To this end, we conducted 100 runs of the \texttt{mies}, \texttt{cma-IH}, as well as the \texttt{cma-wM} (set to box-boundaries of $[-1000,1000]^{64}$) | and investigated the attained optima and their distribution. 
The behavior is consistent across the four instances, and thus we will refer to the specific instance of $\{c=10,E=10\}$, which is depicted in Figure \ref{fig:resIRR10}. \\
We report on our main observations:
\begin{itemize}
    \item As evident in Figure \ref{fig:resIRR10}, there exists a large volume of high-quality solutions which possess objective function values within 0.1\% of the best-attained value (the histogram of the 300 runs is presented on the left), while having diverse underlying decision vectors (the mutual $L_1$ distances among the attained solution vectors are presented by a colormap on the right). 
    The vast majority of the runs converged to unique solution vectors, which are always non-neighboring on the integer lattice (that is, require non-trivial integer variations). Notably, the \texttt{cma-wM} obtained a solution set of the lowest diversity and of the highest quality -- and yet, its solution vectors are all unique. 
    \item CPLEX does not always locate the best optimizers in the original settings. When re-setting the relative MIP optimality gap to $10^{-5}$ (\texttt{cplex.epgap = 1e-5}) it succeeds to locate them.
    \item The \texttt{cma-wM} consistently outperforms the other MIESs with statistical significance. 
\end{itemize}
We consider the above observations as \textbf{numerical evidence for corroborating the hypothesis that the integer RotEllipse+RotEllipse landscape is multimodal.}
%\td{TODO: joint histogram of CPLEX+MIESs \& L1 Distance Matrix} 
\section{Discussion and Summary}\label{sec:summary}
The proven undecidability of unbounded MIQCQP problems constituted the basis of the current study, which given this potential weakness, addressed research questions involving WBO and BBO approaches targeting that class of problems. 
Our empirical study focused on the strengths and limitations of the WBO CPLEX solver versus two BBO MIESs.
We looked at three main factors of problem complexity: the constraints level (``parabolic tightness''), the interaction between decision variables (separability), and the conditioning of the Hessian matrices of the objective and constraint functions. 
These difficulties were represented in a set of nine selected test-cases, which were scaled by the constraint-leveling and conditioning to yield 216 problem instances. 
Next, we summarize the key takeaways of our study:
\begin{enumerate}
\item \label{sum:cplex} The CPLEX solver is significantly challenged by MIQCQP problems with unbounded integer decision variables (reflected by consistent termination in \textit{timeouts}), but usually provides optima of the best quality (even if obtained in \textit{timeouts}). Overall, it is outperformed in only 13\% of the cases.  
The RotEllipse+RotEllipse test-case (TC-4) is an exception to the general trend of \textit{timeouts}. 
We claim that its strongly convex structure, which is rooted in the pair of aligned ellipsoids $\left(\mathbf{H}_f,\mathbf{H}_g\right)$ (see the commentary to Figure \ref{fig:landscapes-gallery}), enables CPLEX to find accurate solutions swiftly.
%It experiences difficulties when the search landscape possesses a \textbf{ridge structure} (Cigar across all conditioning), and on the other hand, swiftly handles landscapes with richer multimodal structures (high-conditioned RotEllipse). 
%\item The easiest instances for the CPLEX in terms of termination criteria (RotE+RotE with $c>10$) are found to be the hardest for the \texttt{mies}.
\item The two MIESs were shown to be well-suited for unbounded integer search tasks, relying on their self-adaptive mutation operators. The comparison between the \texttt{mies} to the \texttt{cma-IH} was unfair in the first place, since the former lacks the ability to produce correlated mutations, which strongly underlies the latter.
Accordingly, the \texttt{mies} could not effectively handle the non-separable landscapes, especially when highly conditioned, unlike the \texttt{cma-IH}, which performed well but was challenged in integer precision when highly conditioned.
\item The \texttt{cma-IH} always completes 10 sequential runs with precise results on the 24 Cigar+Cigar instances, preceding the 1 hour time-limit of the single CPLEX run. 
This empirical observation makes it an attractive first-choice on such a test-case.
\item The combined effect of unboundedness and optima translation leads to a significant decline in the performance of CPLEX, causing it to fail in 85\% of the cases. 
At the same time, the MIESs maintain consistent performance, demonstrating their self-adaptation capabilities.
%The amplification of the unboundedness effect by the optima translation causes CPLEX to fail miserably and exhibit inferior performance on 85\% of the cases.
\item Integer mutations using the double-geometric distribution showed great potential when examining the integer delta measure of solutions obtained by the \texttt{mies}. 
Also, the \texttt{mies} could benefit from a stronger strategy for the continuous search, even on separable landscapes.
%\item The easiest instances for the \texttt{cma-IH} (in terms of termination criteria are challenging for the CPLEX (Cigar+Cigar across all instances); in fact, the \texttt{cma-IH} always completes 10 sequential runs with precise results prior to the 1hour time-limit of a single CPLEX run.
%\item The \texttt{cma-IH} performed at its best whenever the two Hessian matrices coincided, $\mathbf{H}_f=\mathbf{H}_g$. This is to be expected when the \textit{cost function} \eqref{eq:cost} encapsulates two Hessian matrices, even if the penalty term is a squared form of one of them (yielding a \textit{quartic} polynomial of $\vec{x}$ with a quadratic term). That is, convergence is effectively achieved when the mutation operator reaches its full potential upon successfully learning an accurate covariance matrix, which coherently addresses the two components of the cost function. 
\item The numerical observations were consistent when testing the algorithms on the all-integer version of this test-suite (i.e., the 216 equivalent IQCQP problem instances).
\item The hypothesis stating that the integer RotEllipse+RotEllipse landscape is multimodal was numerically validated. 
\item The Hadamard Ellipse presents difficulty when taking part in either the objective or the constraint function, causing CPLEX to fail to find solutions that the \texttt{cma-IH} can identify. This Hessian has been previously reported to be deceptive in different contexts of optimization (\cite{Shir-Yehudayoff_TCS2020} and \cite{moMIQC2024}), likely due to its special form (eigenvalue spectrum of only $\pm 1$) and its distinct impact as a rotation matrix (see also the commentary to Figure \ref{fig:landscapes-gallery}). \\
We hypothesize whether the failure of the CPLEX is rooted in the facts that the integer landscapes are multimodal and that the origin (i.e., the zero vector) constitutes a sub-optimal solution of high quality when the HadEllipse is involved. The CPLEX evidently outputs the origin as its solution in multiple cases.
\end{enumerate}

Evidently, it is not a clear cut decision whether to employ the WBO CPLEX or the BBO MIESs, since their relative performance largely depends on the characteristics of the quadratic forms. 
Moreover, conditioning and separability are not intuitive factors in determining the MIQCQP complexity. 
Ridge-like versus convex structures %(i.e., low- versus high multimodality in this context) 
can pose mirrored degrees of challenge for WBO versus BBO | as was evident in the mirrored performance patterns of CPLEX versus the MIESs, which was summarized in point \ref{sum:cplex} above.

\medskip

Finally, we list possible directions for future research:
\begin{itemize}
    \item Study the impact of constraints handling on the MIESs' performance over MIQCQP problems.  
    \item Investigate the effectiveness of the double-geometric distribution when in concert with a correlation mechanism in order to handle the non-separable landscapes.
    \item Explain any peculiar or unusual observations in the performance patterns when considering constraint level $E$ and conditioning $c$. 
    \item Address MIQCQP from the theoretical perspective and aim to obtain first results on MIESs' step-size behavior and/or convergence on such MI problems.
\end{itemize}

%but it must be viewed together with the structure of the Hessian matrix. 
%\td{Moreover, our study shows that the isotropic mutation of integer variables might be an inferior strategy in case of rotated Hessian matrices, and mechanisms such as the covariance matrix adaptation \cite{ESchapter2018} might have to be adopted to make BBO meta-heuristics 
%On the side of WBO algorithms, the observed problem with loosely bounded decision variables should be regarded as an incentive for further research, but also their under-achieved coverage for non-separable problems deserves attention.}

%\bibliographystyle{ieeetr}
%\bibliography{oshir}

\end{document}